\documentclass[12pt]{article}
\usepackage{latexsym,amsmath,amssymb,amsfonts,epsfig,graphicx,cite,psfrag}
\usepackage{ifpdf}

\newtheorem{thm}{Theorem}[section]
\newtheorem{lem}[thm]{Lemma}
\newtheorem{cor}[thm]{Corollary}

\newtheorem{defn}[thm]{Definition}

\newtheorem{exmp}[thm]{Example}

\def\qed{\hfill \rule{4pt}{7pt}}
\def\pf{\noindent {\bf Proof}:\ }

\newcommand{\color}[2][x]{} %deals with undefined mnacro in pstex_t files

\def\Z{{\mathbb Z}}
\def\C{{\mathbb C}}

\def\mysmallstrut{\hbox{\vrule height5.0pt depth2pt width0pt}}
\def\mysmallbox{\hbox to 7.0pt}
\def\mystrut{\hbox{\vrule height9.0pt depth2pt width0pt}}
\def\mybox{\hbox to 11.0pt}
\def\norulefill{\leaders\hrule height0pt\hfill}
\def\nr#1{\multispan{#1}\norulefill}
\def\hr#1{\multispan{#1}\hrulefill}

\begin{document}

\title{A hive model determination of
multiplicity-free Schur function products
and skew Schur functions}

\author{Donna Q.J. Dou\thanks{email: qjdou@cfc.nankai.edu.cn}
and Robert L. Tang\thanks{email: tangling@cfc.nankai.edu.cn}
\\
Center for Combinatorics, Nankai University,\\
Tianjin 300071, China\\
and\\
Ronald C. King\thanks{email: R.C.King@soton.ac.uk}
\\
School of Mathematics, University of Southampton, \\
Southampton SO17 1BJ, England \\
}

\date{27 March 2008}

%\begin{keyword}
%Hive model, Littlewood-Richardson coefficients, skew products,
%multiplicity-free.
%\end{keyword}

\maketitle

\begin{abstract}
The hive model is a combinatorial device that may be used to
determine Littlewood-Richardson coefficients and study their
properties. It represents an alternative to the use of the
Littlewood-Richardson rule. Here properties of hives are used to
determine all possible multiplicity-free Schur function products and
skew Schur function expansions. This confirms the results of
Stembridge~\cite{stem}, Gutschwager~\cite{guts} and Thomas and
Yong~\cite{tho-yon}, and sheds light on the combinatorial origin
of the conditions for being multiplicity-free, as well as
illustrating some of the key features and power of the hive model.
\end{abstract}

\section{Introduction}
\label{sec-introduction}

Throughout this paper we will adopt the notation and terminology on
Schur functions taken from the standard text by Macdonald~\cite{macd}.
The Schur functions $s_\lambda$,
indexed by partitions $\lambda$, form a $\Z$-basis of the ring of
symmetric functions. This basis is orthonormal, and the
corresponding bilinear form enables one to define skew Schur
functions $s_{\lambda/\mu}$, indexed by pairs of partitions
$\lambda$ and $\mu$ with $\mu\subseteq\lambda$. Each such skew Schur
function $s_{\lambda/\mu}$ can be expressed as a sum of Schur
functions $s_\nu$ with non-negative integer multiplicities given by
the well-known Littlewood-Richardson coefficients,
$c_{\mu\nu}^\lambda$. These same coefficients govern the
decomposition of the product $s_\mu s_\nu$ of two Schur functions as
a sum of Schur functions $s_\lambda$.

In \cite{stem}, Stembridge classified the products of Schur
functions that are multiplicity-free, that is those pairs of Schur
functions for which every coefficient in the Schur function
expansion of their product is $0$ or $1$, as in the following theorem:

\begin{thm}[Stembridge~{\rm\cite{stem}}]
\label{thm-stem}
The Schur function product $s_\mu s_\nu$ is multiplicity-free if and only if
one or more of the following is true:
\begin{description}
\item[P0]~~$\mu$ or $\nu$ is the zero partition $0$;
\item[P1]~~$\mu$ or $\nu$ is a one-line rectangle;
\item[P2]~~$\mu$ is a two-line rectangle and $\nu$ is a fat hook (or
vice versa);
\item[P3]~~$\mu$ is a rectangle and $\nu$ is a near-rectangle (or
vice versa);
\item[P4]~~$\mu$ and $\nu$ are rectangles.
\label{description-stemlist}
\end{description}
Here each partition $\lambda$ is to be identified with the corresponding
Young diagram $F^\lambda$, and is said to be a rectangle if $\lambda=(a^p)$ with
$a>0$ and $p>0$, and a fat hook if $\lambda=(a^pb^q)$ with $a>b>0$ and
$p,q>0$. A rectangle $(a^p)$ is a one-line rectangle if $a=1$ or $p=1$,
and a two-line rectangle if $a,p>1$ with $a=2$ or $p=2$. A fat hook
$(a^pb^q)$ is a near rectangle if any one or more of $a-b$, $b$, $p$ or $q$
is equal to $1$.
\end{thm}

Stembridge~\cite{stem} also gave a corresponding theorem applicable to the
case where the lengths of the partitions satisfy $\ell(\mu),\ell(\nu)<n$
and the Schur function product is restricted to the ring, $\Lambda_n$,
of symmetric functions in a finite number of variables $x_1,x_2,\ldots,x_n$.
To do this he invoked, in the particular case $m=\lambda_1$, the use of
the partition $\lambda^*=(m-\lambda_n,m-\lambda_{n-1},\ldots,m-\lambda_1)$,
which may be said to be $m^n$-complementary to $\lambda$.

More recently, Thomas and Yong~\cite{tho-yon} established a
multiplicity-free result for the product of two Schubert classes
$\sigma_\mu\,\sigma_\nu$ in the cohomolgy ring
$H^*(Gr(m,\C^{m+n}),\Z)$ of the Grasmannian $Gr(m,\C^{m+n})$ of
$m$-dimensional subspaces in $\C^{m+n}$. The coefficients in this
product of Schubert classes are just the usual Littlewood-Richardson
coefficients, $c_{\mu\nu}^\lambda$, but this time restricted to the
case $\lambda\subseteq m^n$. This allows, Thomas and Yong's result
to be recast in terms of skew Schur functions. When this is done it
coincides with the multiplicity-free result for skew Schur functions
derived independently by Gutschwager~\cite{guts}.

Here it is convenient to state their common result just for {\em basic}
skew Schur functions, that is
for those cases $s_{\lambda/\mu}$ where the skew Young diagram
$F^{\lambda/\mu}$ has neither empty rows nor empty columns. It
need not be connected. As will be seen, every
skew Schur function is equal, in a rather trivial way, to some basic
skew Schur function. In stating the theorem
it is also convenient to follow Thomas
and Yong in letting the $m^n$-shortness of a partition $\lambda\subseteq m^n$
be the length of the shortest straight line segment of the path
of length $m+n$ from the southwest to northeast corner of $F^{m^n}$
that separates $F^\lambda$ from the $\pi$-rotation of $F^{\lambda^*}$.
With this definition, the result, jointly attributable both to
Gutschwager and to Thomas and Yong, takes the form:

\begin{thm}[Gutschwager~{\rm\cite{guts}}, Thomas and Yong~{\rm\cite{tho-yon}}]
\label{thm-main} The basic skew Schur function $s_{\lambda/\mu}$ is
multiplicity-free if and only if one or more of the following is
true:
\begin{description}
\item[R0]~~$\mu$ or $\lambda^*$ is the zero partition $0$;
\item[R1]~~$\mu$ or $\lambda^*$ is a rectangle of $m^n$-shortness $1$;
\item[R2]~~$\mu$ is a rectangle of $m^n$-shortness $2$ and
$\lambda^*$ is a fat hook (or vice versa);
\item[R3]~~$\mu$ is a rectangle and $\lambda^*$ is a fat hook of $m^n$-shortness $1$
(or vice versa);
\item[R4]~~$\mu$ and $\lambda^*$ are rectangles.
\label{description-gtylist}
\end{description}
where $\lambda^*$ is the $m^n$-complement of $\lambda$ with $m=\lambda_1$
and $n=\lambda'_1$.
\end{thm}

In deriving this result both Gutschwager~\cite{guts}
and Thomas and Yong~\cite{tho-yon}
used the traditional model of Littlewood-Richardson
coefficients~\cite{litt-richa,litt}:
\begin{thm}[Littlewood and Richardson~{\rm\cite{litt-richa}}]
The Littlewood-Richardson coefficient $c^{\lambda}_{\mu\nu}$ is
equal to the number of semistandard Young tableaux of shape
$\lambda/\mu$ and content $\nu$ whose reverse reading word is a
lattice permutation. \label{thm-LRrule}
\end{thm}

Here, our aim is to rederive both Theorems~\ref{thm-stem}
and~\ref{thm-main} by means of a different combinatorial model,
namely the hive model~\cite{kus,buch} for these coefficients:

\begin{thm}[Knutson and Tao~{\rm\cite{kus}}]
The Littlewood-Richardson coefficient $c^{\lambda}_{\mu\nu}$ is
equal to the number of distinct integer LR-hives with boundary edge
labels specified by the partitions $\lambda$, $\mu$ and $\nu$.
\label{thm-LRhives}
\end{thm}

The hive model introduced in~\cite{kus} arose as a reformulation
of a convex polytope model~\cite{bzel}. It was described in more
detail in~\cite{buch}, with an appendix providing a rather simple
bijection between the tableaux of Theorem~\ref{thm-LRrule}
and the hives of Theorem~\ref{thm-LRhives}. For more information about
the hive model, see also~\cite{ktw,ktt}.

There are a number of advantages to the hive model approach, including
the fact that it allows a direct proof that all the
cases enumerated in both Theorem~\ref{thm-stem} and
Theorem~\ref{thm-main} are indeed multiplicity-free. It also lends
itself well to the uniform statement of Theorem~\ref{thm-main},
simultaneously covering both connected and disconnected cases.
What is avoided in the hive model proof that all the indicated
multiplicity-free cases are indeed multiplicity-free
is any recourse to
the non-trivial order filter, introduced by Stembridge and
generalised by Gutschwager, that underlies two families of
inequalities of the form $c^\lambda_{\mu\nu}\geq
c^\rho_{\sigma\tau}$ that are also heavily used by
Thomas and Yong in the form of what they call Stembridge
demolitions. Although we cannot avoid the use of such
inequalities in dealing with all possible non-multiplicity-free
cases, the hive model does allow us to be completely explicit about
the route from the most general cases to those for which a
multiplicity of at least two occurs. In doing so the hive model
offers some insight into the origin of the breakdown of
multiplicity-freeness for both products of Schur functions and
expansions of skew Schur functions. This lies in the fact that
within the appropriate LR-hives there always exists an elementary
hexagon whose interior edge labels are not fixed. The precise nature
of the somewhat numerous conditions for the breakdown of
multiplicity-freeness can then be exposed in each case through
consideration of a single hive diagram.

In the next section, we make some necessary definitions regarding
partitions, skew Schur functions and Littlewood-Richardson
coefficients. In Section~\ref{sec-hive-model} we define integer
$n$-hives and the LR-hives whose enumeration for fixed boundary
labels provides a model for evaluating Littlewood-Richardson
coefficients. A sequence of lemmas regarding LR-hives and subhives
are derived in Section~\ref{sec-LRhives}. These are used in
Sections~\ref{sec-Stembridge} and~\ref{sec-main-result},
respectively, to prove that all the Schur function products and skew
Schur functions listed in Theorems~\ref{thm-stem}
and~\ref{thm-main} are indeed multiplicty-free. The question of
completeness of these lists is tackled in
Sections~\ref{sec-Stembridge-completeness}
and~\ref{sec-main-completeness}, thereby completing the proof of
both theorems. Some final remarks, including a corollary regarding
multiplicity-free products of skew Schur functions, are offered in
Section~\ref{sec-final-remarks}.

\section{Skew Schur functions and Littlewood-Richardson coefficients}
\label{sec-LR}

Let $n$ be a fixed positive integer, and let
$\lambda=(\lambda_1,\lambda_2,\ldots,\lambda_n)$ be a partition of
weight $|\lambda|=\lambda_1+\lambda_2+\cdots+\lambda_n$ and length
$\ell(\lambda)\leq n$. The parts of $\lambda$ are non-negative
integers such that $\lambda_1\geq \lambda_2\geq\cdots\geq\lambda_n$
with $\lambda_i>0$ for all $i\leq\ell(\lambda)$ and $\lambda_i=0$
for all $i>\ell(\lambda)$. Such a partition $\lambda$ specifies a
Young diagram $F^\lambda$ consisting of $|\lambda|$ boxes whose row
lengths are the parts $\lambda_i$ of $\lambda$ and whose column
lengths are the parts $\lambda'_j$ of the conjugate partition
$\lambda'$. Schematically, we have:
$$
F^\lambda=\quad
{\vcenter
 {\offinterlineskip
 \halign{&\mystrut\vrule#&\mybox{\hss$\scriptstyle#$\hss}\cr
  \hr{11}\cr
                      &\lambda_1&\omit& &\omit& &\omit& &\omit& &\cr    %+5
  \hr{11}\cr
                      &\lambda_2&\omit& &\omit& &\omit& &\cr       %+4
  \hr{9}\cr
                      &\lambda_3&\omit& &\omit& &\omit& &\cr       %+4
  \hr{9}\cr
                      &\lambda_4&\omit& &\omit& &\cr          %+3
  \hr{7}\cr
 }}}
\quad=\quad
{\vcenter
 {\offinterlineskip
 \halign{&\mystrut\vrule#&\mybox{\hss$\scriptstyle#$\hss}\cr
  \hr{11}\cr
   &\lambda_1^\prime&&\lambda_2^\prime&&\lambda_3^\prime&&\lambda_4^\prime&&\lambda_5^\prime&\cr    %+5
  \hr{1}&\nr{7}&\hr{3}\cr
                                         & && && && &\cr    %+4
  \hr{1}&\nr{7}&\hr{1}\cr
                                         & && && && &\cr    %+4
  \hr{1}&\nr{5}&\hr{3}\cr
                                            & && && &\cr    %+3
  \hr{7}\cr
 }}}
\qquad\qquad
F^{\lambda'}=\quad
{\vcenter
 {\offinterlineskip
 \halign{&\mystrut\vrule#&\mybox{\hss$\scriptstyle#$\hss}\cr
  \hr{9}\cr
                      &\lambda'_1&\omit& &\omit& &\omit& &\cr    %+4
  \hr{9}\cr
                      &\lambda'_2&\omit& &\omit& &\omit& &\cr       %+4
  \hr{9}\cr
                      &\lambda'_3&\omit& &\omit& &\omit& &\cr       %+4
  \hr{9}\cr
                      &\lambda'_4&\omit& &\omit& &\cr          %+3
  \hr{7}\cr
                      &\lambda'_5&\cr          %+1
  \hr{3}\cr
 }}}
\quad=\quad
{\vcenter
 {\offinterlineskip
 \halign{&\mystrut\vrule#&\mybox{\hss$\scriptstyle#$\hss}\cr
  \hr{9}\cr
                      &\lambda_1&&\lambda_2&&\lambda_2&&\lambda_4&\cr    %+4
  \hr{1}&\nr{7}&\hr{1}\cr
                      & && && && &\cr       %+4
  \hr{1}&\nr{7}&\hr{1}\cr
                      & && && && &\cr       %+4
  \hr{1}&\nr{5}&\hr{3}\cr
                      & && && &\cr          %+3
  \hr{1}&\nr{1}&\hr{5}\cr
                      & &\cr          %+1
  \hr{3}\cr
 }}}
$$

It is sometimes convenient to write
$\lambda=(\lambda_1,\lambda_2,\ldots,\lambda_n)$ in the form
$\lambda=(a^p,b^q,\ldots)$ if $\lambda_i=a$ for $1\leq i\leq p$,
$\lambda_i=b$ for $p<i\leq p+q$, etc., with $a>b>\cdots>0$ and
$p,q,\ldots>0$. In addition, for any pair of partitions $\lambda$
and $\mu$ we define $\lambda+\mu$ to be the partition obtained by
adding corresponding parts of $\lambda$ and $\mu$, and
$\lambda\cup\mu$ is the partition obtained by arranging all the
parts of $\lambda$ and $\mu$ in weakly decreasing order.

We write $\mu\subseteq\lambda$ if all the boxes of $F^\mu$ are
contained in $F^\lambda$, that is to
say $\mu_i\leq\lambda_i$ for all $i$, or equivalently,
$\mu'_j\leq\lambda'_j$ for all $j$. In such a case
the corresponding skew diagram $F^{\lambda/\mu}$ is the diagram
obtained by deleting from $F^\lambda$ all the boxes
of $F^\mu$. Schematically, we have:
$$
F^{\lambda}=\quad {\vcenter
 {\offinterlineskip
 \halign{&\mystrut\vrule#&\mybox{\hss$\scriptstyle#$\hss}\cr
  \hr{23}\cr
  &~~\mu_1&\omit& &\omit& &\omit& &\omit& &&~~~~~\lambda_1-\mu_1&\omit& &\omit& &\omit& &\omit& &\omit& &\cr
  \hr{23}\cr
  &~~\mu_2&\omit& &\omit& &\omit& &&~~~~~\lambda_2-\mu_2&\omit& &\omit& &\omit& &\cr
  \hr{17}\cr
  &~~\mu_3&\omit& &\omit& &&~~~~~\lambda_3-\mu_3&\omit& &\omit& &\omit& &\omit& &\cr
  \hr{17}\cr
  &~~\mu_4&\omit& &\omit& &&~~~~~\lambda_4-\mu_4&\omit& &\omit& &\cr
  \hr{13}\cr
   }}}
\qquad F^{\lambda/\mu}=\quad {\vcenter
 {\offinterlineskip
 \halign{&\mystrut\vrule#&\mybox{\hss$\scriptstyle#$\hss}\cr
  \nr{4}&\hr{13}\cr
  \omit& &\omit& &&~~~~~\lambda_1-\mu_1&\omit& &\omit& &\omit& &\omit& &\omit& &\cr
  \nr{2}&\hr{15}\cr
  \omit& &&~~~~~\lambda_2-\mu_2&\omit& &\omit& &\omit& &\cr
  \hr{11}\cr
  &~~~~~\lambda_3-\mu_3&\omit& &\omit& &\omit& &\omit& &\cr
  \hr{11}\cr
  &~~~~~\lambda_4-\mu_4&\omit& &\omit& &\cr
  \hr{7}\cr
   }}}
$$

Just as to each partition $\lambda$ there corresponds a Schur function
$s_\lambda$, so to each pair of partitions $\lambda$ and $\mu$
with $\mu\subseteq\lambda$ there corresponds a skew Schur function
$s_{\lambda/\mu}$~\cite{litt,macd}.
This may be defined by noting first that there exists
a symmetric bilinear form $\langle\cdot\,,\,\cdot\rangle$ on the ring
of symmetric functions $\Lambda$
such that $\langle s_\mu\,,\,s_\nu\rangle=\delta_{\mu\nu}$,
and then defining $s_{\lambda/\mu}$ by the relations~\cite{macd}
\begin{equation}
\langle s_{\lambda/\mu}\,,\,s_\nu\rangle=\langle s_\lambda\,,\,s_\mu\,s_\nu\rangle\,,
\label{eq-skew-sfn}
\end{equation}
for all partitions $\nu$.

Since the Littlewood-Richardson coefficients
$c_{\mu\nu}^\lambda$ arise as the multiplicities in
the expansion of the Schur function product
\begin{equation}
s_\mu\, s_\nu = \sum_{\lambda}\ c^{\lambda}_{\mu\nu}\,s_{\lambda}\,,
\label{eq-sfn-LRcoeff}
\end{equation}
it follows from (\ref{eq-skew-sfn})
that they must also arise as the multiplicities in the skew Schur function expansion
\begin{equation}
s_{\lambda/\mu} = \sum_{\nu}\ c^{\lambda}_{\mu\nu}\,s_{\nu}\,.
\label{eq-skew-sfn-LRcoeff}
\end{equation}

The Littlewood-Richardson rule implies that $c^\lambda_{\mu\nu}$
can only be non-zero if
\begin{equation}
|\lambda|=|\mu|+|\nu|
\qquad\hbox{and}\qquad
\ell(\mu),\ell(\nu)\leq \ell(\lambda)\leq\ell(\mu)+\ell(\nu).
\label{eq-weight-length}
\end{equation}

Although it is by no means obvious from the Littlewood-Richardson rule,
the Littlewood-Richardson coefficients satisfy a number of
symmetry properties, including:
\begin{equation}
  c_{\mu\nu}^\lambda=c_{\nu\mu}^\lambda
  \qquad\hbox{and}\qquad
   c_{\mu'\nu'}^{\lambda'}=c_{\mu\nu}^\lambda\,.
\label{eq-symmetry-LRcoeff}
\end{equation}

Moreover,
for all partitions $\lambda$, $\mu$ and $\nu$
and all non-negative integers $a$, $b$ and $c$ with $a=b+c$, we have
\begin{equation}
c^{\lambda+(1^a)}_{\mu+(1^b),\nu+(1^c)}\geq
  c^\lambda_{\mu\nu}
  \qquad\hbox{and}\qquad
c^{\lambda\cup(a)}_{\mu\cup(b),\nu\cup(c)}\geq
  c^\lambda_{\mu\nu}
  \,.
\label{eq-add-col-row}
\end{equation}
These inequalities, which are related by conjugacy, have been
derived in~\cite{guts}, as a generalisation of the $c=0$ case given
in~\cite{stem}. Although a direct proof of both inequalities may be
based on the hive model, we do not present the proof here.

It is useful to note~\cite{sta} that
\begin{equation}
s_\lambda=s_{\lambda^\pi}
\qquad\hbox{and}\qquad
s_{\lambda/\mu}=s_{(\lambda/\mu)^\pi} \,,
\label{eq-skew-rotate}
\end{equation}
where $F^{\lambda^\pi}$ and $F^{(\lambda/\mu)^\pi}$ are obtained by
rotating $F^\lambda$ and $F^{\lambda/\mu}$, respectively, through
$\pi$ radians.

\begin{exmp}
If $\lambda=(432)$ and $\mu=(2)$, then the $\pi$-rotations of $F^\lambda$ and $F^{\lambda/\mu}$
take the form:
$$F^{\lambda}=\
{\vcenter
 {\offinterlineskip
 \halign{&\mysmallstrut\vrule#&\mysmallbox{\hss$#$\hss}\cr
  \hr{9}\cr
  & && && && &\cr
  \hr{9}\cr
  & && && &\cr
  \hr{7}\cr
  & && &\cr
  \hr{5}\cr
   }}}
\Rightarrow%\qquad
F^{\lambda^\pi}=\
 {\vcenter
 {\offinterlineskip
 \halign{&\mysmallstrut\vrule#&\mysmallbox{\hss$#$\hss}\cr
  \nr{4}&\hr{5}\cr
  \omit&*&\omit&*&& && &\cr
  \nr{2}&\hr{7}\cr
  \omit&*&& && && &\cr
  \hr{9}\cr
  & && && && &\cr
  \hr{9}\cr
   }}}
\quad\hbox{and}\quad
F^{\lambda/\mu}=\
{\vcenter
 {\offinterlineskip
 \halign{&\mysmallstrut\vrule#&\mysmallbox{\hss$#$\hss}\cr
  \nr{4}&\hr{5}\cr
  \omit&*&\omit&*&& && &\cr
  \hr{9}\cr
  & && && &\cr
  \hr{7}\cr
  & && &\cr
  \hr{5}\cr
   }}}
\Rightarrow%\qquad
F^{(\lambda/\mu)^\pi}=\
 {\vcenter
 {\offinterlineskip
 \halign{&\mysmallstrut\vrule#&\mysmallbox{\hss$#$\hss}\cr
  \nr{4}&\hr{5}\cr
  \omit&*&\omit&*&& && &\cr
  \nr{2}&\hr{7}\cr
  \omit&*&& && && &\cr
  \hr{9}\cr
  & && &\cr
  \hr{5}\cr
   }}}
$$
\noindent so that we have $s_{432}=s_{444/21}$ and
$s_{432/2}=s_{442/21}$. \label{exmp-rot}
\end{exmp}

As far as $m^n$-complements are concerned
\begin{equation}
s_{m^n/\lambda}=\left\{ \begin{array}{ll}
  s_{\lambda^*}&~~\hbox{if}~~\lambda\subseteq m^n\,;\cr
  0&~~\hbox{otherwise}\,,\cr
  \end{array} \right.
\label{eq-skew-complement}
\end{equation}
where $\lambda^*_k=m-\lambda_{n-k+1}$ for $k=1,2,\ldots,n$.

An important consequence of this rather trivial observation
is that for $\mu\subseteq\lambda\subseteq m^n$ we have
\begin{equation}
\begin{array}{lclclcl}
 (s_{\lambda/\mu}\,,\,s_\nu)&=&(s_{\lambda}\,,\,s_\mu\,s_\nu)
 &=&(s_{m^n/\lambda^*}\,,\,s_\mu\,s_\nu)\cr
 &=&(s_{m^n}\,,\,s_{\lambda^*}\,s_\mu\,s_\nu)
 &=&(s_{m^n/\nu}\,,\,s_{\lambda^*}\,s_\mu)
 &=&(s_{\nu^*}\,,\,s_{\lambda^*}\,s_\mu)\cr
\end{array}
\label{eq-skew-schubert}
\end{equation}
with the result non-zero only if $\nu\subseteq m^n$.
It is this identity which enables us to conclude that
the skew Schur function $s_{\lambda/\mu}$ is
multiplicity-free if and only if the product of Schubert
classes $\sigma_{\lambda^*}\,\sigma_\mu$ is
multiplicity-free.

\begin{exmp}
By way of illustration, if $m=9$, $n=5$, $\lambda=(99666)$ and
$\mu=(552)$ then $\lambda^*=(333)$ and $F^{\lambda/\mu}$
and its $m^n$-complement take the form:
$$F^{\lambda/\mu}=\
{\vcenter
 {\offinterlineskip
 \halign{&\mysmallstrut\vrule#&\mysmallbox{\hss$#$\hss}\cr
  \nr{10}&\hr{9}\cr
  \omit& &\omit& &\omit& &\omit& &\omit& && && && && &\cr
  \nr{10}&\hr{9}\cr
  \omit& &\omit& &\omit& &\omit& &\omit& && && && && &\cr
  \nr{4}&\hr{15}\cr
  \omit& &\omit& && && && && &\cr
  \hr{13}\cr
  & && && && && && &\cr
  \hr{13}\cr
  & && && && && && &\cr
  \hr{13}\cr
   }}}
\qquad\qquad
F^{\mu}\cup F^{\lambda^{*\pi}}=\
 {\vcenter
 {\offinterlineskip
 \halign{&\mysmallstrut\vrule#&\mysmallbox{\hss$#$\hss}\cr
  \hr{19}\cr
  & && && && && && &\omit& &\omit& &\omit& &\cr
  \hr{11}&\nr{7}&\hr{1}\cr
  & && && && && && &\omit& &\omit& &\omit& &\cr
  \hr{11}&\nr{1}&\hr{7}\cr
  & && && &\omit& &\omit& &\omit& && && && &\cr
  \hr{5}&\nr{7}&\hr{7}\cr
  & &\omit& &\omit& &\omit& &\omit& &\omit& && && && &\cr
  \hr{1}&\nr{11}&\hr{7}\cr
  & &\omit& &\omit& &\omit& &\omit& &\omit& && && && &\cr
  \hr{19}\cr
   }}}
$$
\end{exmp}
In this example it can be seen that the sequence of straight line
segments of the two paths from the southwest to northeast corners of
the $m^n=9^5$ rectangle that border the inner and outer boundaries
of $F^{\lambda/\mu}$ are given by $s_{in}=(2,2,1,3,2,4)$ and
$s_{out}(6,3,3,2)$, respectively, where the terminology
of~\cite{guts} has been adopted. The shortness of $\mu$ and
$\lambda^*$ as defined in~\cite{tho-yon} are just the smallest
components, $1$ and $2$ respectively, of these two sequences. 

A futher useful fact about skew Schur functions is that
\begin{equation}
s_{\lambda/\mu}=s_{\hat\lambda/\hat\mu}\,,
\label{eq-skew-basic}
\end{equation}
where $F^{\hat\lambda/\hat\mu}$ is the skew Young diagram
obtained from $F^{\lambda/\mu}$ by deleting any empty rows,
that is those for which $\lambda_i=\mu_i$,
and any empty columns, that is those for which $\lambda'_j=\mu'_j$.
The skew Schur function $s_{\hat\lambda/\hat\mu}$ is said to be
{\it basic}. This identity therefore allows each skew Schur function
to be expressed as a basic skew Schur function. It should be noted
that if $s_{\lambda/\mu}$ is itself basic, then
$\mu_i<\lambda_i$ and $\mu_i\leq\lambda_{i+1}$ for
$i=1,2,\ldots,\ell(\lambda)-1$, with $\ell(\mu)<\ell(\lambda)$.
If we just have $\mu_i<\lambda_i$ for all $i=1,2,\ldots,\ell(\lambda)$,
then we say that $s_{\lambda/\mu}$ is {\it row-basic}.

\begin{exmp}
\label{exmp-basic}
In the case $\lambda=(985333)$ and $\mu=(755321)$ the construction
of $F^{\hat\lambda/\hat\mu}$ from $F^{\lambda/\mu}$ is illustrated
by:
$$
F^{\lambda/\mu}=\quad
{\vcenter
 {\offinterlineskip
 \halign{&\mysmallstrut\vrule#&\mysmallbox{\hss$#$\hss}\cr
  \nr{14}&\hr{5}\cr
  \omit&\cdot&\omit&*&\omit&*&\omit&\cdot&\omit&\cdot&\omit&*&\omit&*&& && &\cr
  \nr{10}&\hr{9}\cr
  \omit&\cdot&\omit&*&\omit&*&\omit&\cdot&\omit&\cdot&& && && &\cr
  \nr{10}&\hr{7}\cr
  \omit&\cdot&\omit&\cdot&\omit&\cdot&\omit&\cdot&\omit&\cdot&\cr
  \nr{6}&\hr{5}\cr
  \omit&\cdot&\omit&\cdot&\omit&\cdot&\cr
  \nr{4}&\hr{3}\cr
  \omit&\cdot&\omit&*&& &\cr
  \nr{2}&\hr{5}\cr
  \omit&\cdot&& && &\cr
  \nr{2}&\hr{5}\cr
%  \omit&\cdot&& &\cr
%  \hr{5}\cr
   }}}
\qquad\Rightarrow\qquad
  F^{\hat\lambda/\hat\mu}=\quad
 {\vcenter
 {\offinterlineskip
 \halign{&\mysmallstrut\vrule#&\mysmallbox{\hss$#$\hss}\cr
  \nr{8}&\hr{5}\cr
  \omit&*&\omit&*&\omit&*&\omit&*&& && &\cr
  \nr{4}&\hr{9}\cr
  \omit&*&\omit&*&& && && &\cr
  \nr{2}&\hr{9}\cr
  \omit&*&& &\cr
  \hr{5}\cr
  & && &\cr
  \hr{5}\cr
 % & &\cr
 % \hr{3}\cr
   }}}
$$
\noindent Hence $\hat\lambda=(6522)$ and $\hat\mu=(421)$, so that
$s_{985333/755321}=s_{6522/421}$, with $s_{6522/421}$ basic, but in
this instance not connected.
\end{exmp}

If $F^{\lambda/\mu}$ is not connected, but consists of two
components $F^{\theta}$ and $F^{\phi}$ that have no edge in common,
and may themselves be either Young diagrams or skew {\it Young}
diagrams, then
\begin{equation}
s_{\lambda/\mu}=s_{\theta}\,s_{\phi}.
\label{eq-disconnected}
\end{equation}

\begin{exmp}
\label{exmp-disconnected}
In the case $\lambda/\mu=65221/421$ we have the disconnected diagram
$$
F^{\lambda/\mu}=\quad
{\vcenter
 {\offinterlineskip
 \halign{&\mysmallstrut\vrule#&\mysmallbox{\hss$#$\hss}\cr
  \nr{8}&\hr{5}\cr
  \omit&*&\omit&*&\omit&*&\omit&*&& && &\cr
  \nr{4}&\hr{9}\cr
  \omit&*&\omit&*&& && && &\cr
  \nr{2}&\hr{9}\cr
  \omit&*&& &\cr
  \hr{5}\cr
  & && &\cr
  \hr{5}\cr
  & &\cr
  \hr{3}\cr
   }}}
\qquad\Rightarrow\qquad
F^{\theta}=\quad
{\vcenter
 {\offinterlineskip
 \halign{&\mysmallstrut\vrule#&\mysmallbox{\hss$#$\hss}\cr
  \nr{2}&\hr{3}\cr
  \omit&*&& &\cr
  \hr{5}\cr
  & && &\cr
  \hr{5}\cr
  & &\cr
  \hr{3}\cr
   }}}
\qquad
F^{\phi}=\quad
{\vcenter
 {\offinterlineskip
 \halign{&\mysmallstrut\vrule#&\mysmallbox{\hss$#$\hss}\cr
  \nr{4}&\hr{5}\cr
  \omit&*&\omit&*&& && &\cr
  \hr{9}\cr
  & && && &\cr
  \hr{7}\cr
   }}}
$$
\noindent from which it can be seen that $\theta=221/1$ and $\phi=43/2$. Hence
$s_{65221/421}=s_{221/1}\, s_{43/2}$.
\end{exmp}

\section{The hive model}
\label{sec-hive-model}

An $n$-hive is an array of numbers $a_{ij}$, with $0\leq i,\,j,\,
i+j\leq n$, placed at the vertices of an equilateral triangular
graph. Typically, for $n=4$ their arrangement is as shown below.
Such an $n$-hive is said to be an integer hive if all of its entries
are non-negative integers.
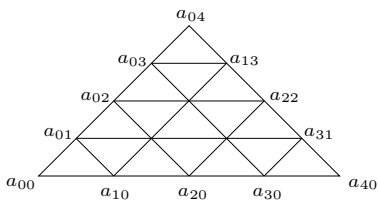
\begin{center}
\input{1-1.TpX}
\end{center}

Neighbouring entries define two distinct types of triangles and
neighbouring triangles define three distinct types of rhombus:

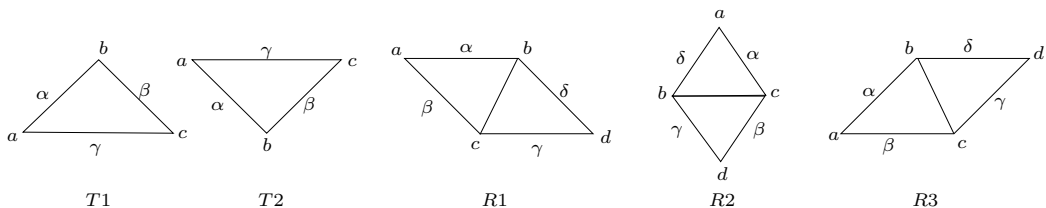
\vspace{-5mm}
\begin{center}
\input{1-2.TpX}
\label{eq-rhombi}
\end{center}
\vspace{-5mm}
each with its own constraint conditions.

In each rhombus, with the labelling as shown above, the hive
condition takes the form:
\begin{equation}
b+c\geq a+d.
\label{eq-hive-condition-1}
\end{equation}

In what follows we make use of edge labels more often than vertex
labels. Each edge in the hive is labelled by means of the
difference, $\epsilon=q-p$, between the labels, $p$ and $q$, of the
two vertices connected by this edge, with $q$ always to the right of
$p$. Thus in both triangles $T1$ and $T2$, we have $\alpha=b-a$,
$\beta=c-b$ and $\gamma=c-a$, so that in each case
\begin{equation}
\label{eq-triangle-condition}
  \alpha+\beta=\gamma\,.
\end{equation}
Similarly, in the case of all three of the above rhombi, $R1$, $R2$ and
$R3$, we have
\begin{equation}
\alpha+\delta=\beta+\gamma\,,
\label{eq-rhombus-condition}
\end{equation}
and the hive conditions
take the form:
\begin{equation}
\alpha\geq \gamma \quad \mbox{and} \quad \beta \geq \delta\,,
\label{eq-hive-condition-2}
\end{equation}
where, of course, either one of these conditions implies the other,
and
\begin{equation}
\alpha=\gamma \quad \mbox{if and only if } \quad \beta=\delta\,.
\label{eq-rhombus-equal-edge}
\end{equation}

We are now in a position to define LR-hives:

\begin{defn}
Let $n$ be a positive integer, and let $\lambda, \mu$ and $\nu$
be any partitions for which $\ell(\lambda),\ell(\mu),\ell(\nu)\leq n$ and
$|\mu|+|\nu|=|\lambda|$. An LR-hive is any integer $n$-hive
with its vertex labels satisfying the hive conditions
\eqref{eq-hive-condition-1} and its boundary vertex labels
given by $a_{00}=0$,
$a_{0,i}=\nu_1+\cdots+\nu_i$, $a_{j,n-j}=|\nu|+\mu_1+\cdots+\mu_j$ and
$a_{k,0}=\lambda_1+\cdots+\lambda_k$ for $i,j,k=1,2,\cdots,n$.

Equivalently, its edge labels satisfy \eqref{eq-triangle-condition},
\eqref{eq-rhombus-condition} and \eqref{eq-hive-condition-2} for all
constituent triangles of type $T1$ and $T2$, and rhombi of type
$R1$, $R2$ and $R3$, and its boundary edge labels are given by
$\lambda_i$, $\mu_j$ and $\nu_k$ for $i,j,k=1,2\ldots,n$.
Schematically, we have: \label{defn-LRhives}
\end{defn}

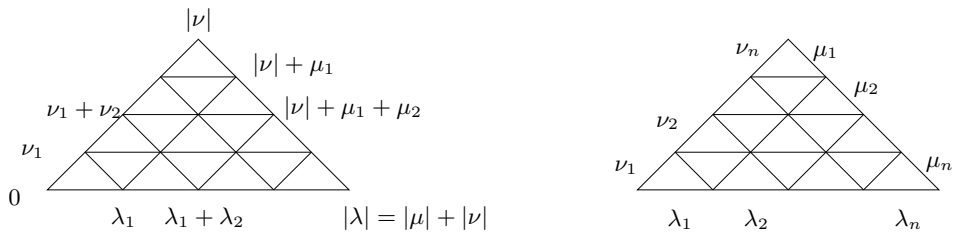
\begin{center}
\input{1-3.TpX}
\end{center}
\vspace{-3mm}

The labelling has been given first in terms of vertex labels and
then in terms of edge labels.
The right hand edge labelling scheme is the one that we will adopt for all
subsequent LR-hives, and it is precisely this type of hive whose
enumeration determines the Littlewood-Richardson coefficient $c_{\mu\nu}^\lambda$
as in Theorem~\ref{thm-LRhives}.

\section{Some properties of LR-hives and subhives}
\label{sec-LRhives}

All edge labels in any integer hive are, of course, integers. In the
case of an LR-hive, these integer labels are necessarily
non-negative and along any straight line parallel to a boundary they
weakly decrease in one particular direction. This is clearly true on
the boundary, since the labels are all parts of partitions. It is
also true of all interior edge labels as may be seen from the
following pair of diagrams of an arbitrary LR-hive:

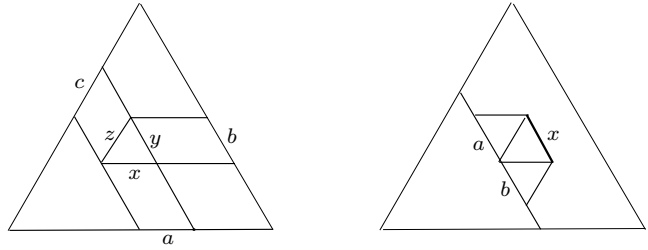
\vspace{-5mm}
\begin{center}
\input{mpf-interior-edges-new.TpX}
\end{center}
\vspace{-3mm}

\noindent
In the first of these diagrams
the hive conditions \eqref{eq-hive-condition-2} applied to the rhombi
of type $R1$, $R1$, $R2$ that constitute the corridors between the
interior edges with labels $x$, $y$, $z$ and the
boundary edges with labels $a$, $b$, $c$, respectively,
imply that $x\geq a\geq0$, $y\geq b\geq0$, $z\geq c\geq0$.
Thus all interior edge labels are non-negative,
as required. In the second diagram, the same hive
conditions applied to rhombi of type $R1$ and $R2$ imply that
$a\geq x$ and $x\geq b$, so that $a\geq b$. This weakly
decreasing condition may readily be extended to cover
all edge labels on the straight line containing $a$ and $b$.
Analogous results apply to edge labels along any straight line
parallel to one or other of the three hive boundaries.

We now establish some properties of subhives of any given
LR-hive. These properties will play an important role in the proof
of our two main theorems.

\begin{lem}
\label{lem-edge-labels}
Each of the following diagrams represents a subhive of an LR-hive.
The subhive may be oriented in any manner within the full hive.
If each edge signified by a solid line is assigned some
fixed label, then these labels are sufficient to determine the labels of
all the remaining edges signified by dashed lines.
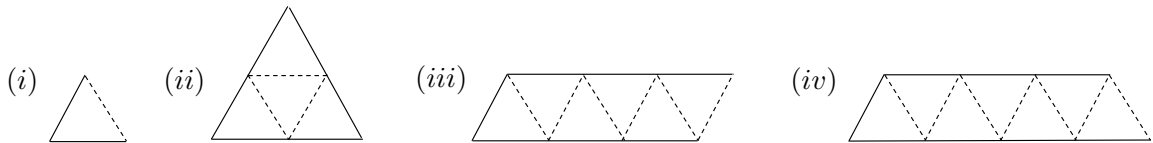
\begin{center}
\input{mpf-iv-lemma.TpX}
\end{center}
\end{lem}

\pf
In each case the repeated use of the triangle condition~\eqref{eq-triangle-condition}
is sufficient to fix all the unassigned edge labels. Case (i) is covered by the fact
that~\eqref{eq-triangle-condition} fixes any one edge label of an
elementary triangle in terms of the other two. Applying this to each of the
the three elementary subtriangles of (ii) then fixes the labels of the three
interior edges. In the cases (iii) and (iv) one can successively determine the
labels on all the dashed line edges by the application of~\eqref{eq-triangle-condition}
to each elementary triangle taken in turn from left to right along each of
these two diagrams.
\qed

\begin{lem}
\label{lem-equal-edge-labels}
Each of the following diagrams represents a subhive of an LR-hive,
which for illustrative purposes has been given a specific orientation.
If each edge signified by a solid line is assigned the
labels $a,b,c,\ldots$ as indicated, then these labels are sufficient
to fix the labels of all the remaining edges signified by dashed lines,
including those dashed line edge labels that have been indicated.
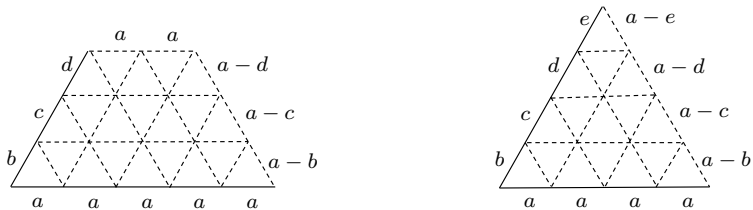
\begin{center}
\input{mpf-equal-edge-lemma.TpX}
\end{center}
\end{lem}

\pf Here it is helpful to consider the following subhive:
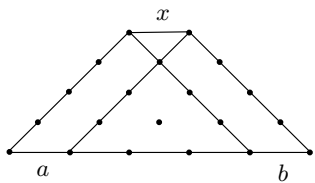
\vspace{-5mm}
\begin{center}
\input{mpf-axa.TpX}
\end{center}
\vspace{-5mm}

\noindent Within this diagram the repeated application of the
rhombus conditions~\eqref{eq-hive-condition-2} to the sequence of
rhombi of type $R3$ between $a$ and $x$, and then a sequence of type
$R1$ between $x$ and $b$, shows that $a\geq x\geq b$. It follows
that if $a=b$ then $x=a$. Applying this to the isosceles trapezium
and the equilateral triangle shows that all the horizontal edges of
both diagrams must have label $a$. Then each of the subhives is seen
to consist of a sequence of thin strip subhives of the case (iv),
dealt with in Lemma~\ref{lem-edge-labels}. It follows that all the
edge labels are fixed not only in each thin strip, but also in both
the isosceles trapezium and the equilateral triangle. With the
orientation as shown, the triangle
condition~\eqref{eq-triangle-condition} and equal edge rhombus
condition~\eqref{eq-rhombus-equal-edge} imply that if the edge
labels on the left hand boundaries are given by $b,c,\ldots$ then
those on the right hand boundary must be $a-b,a-c,\ldots$. For
different orientations $a-b,a-c,\ldots$ must be replaced by either
$a+b,a+c,\ldots$ or $b-a,c-a,\ldots$.

\begin{lem}
\label{lem-hexagon} In each of the following diagrams the elementary
hexagon represents a subhive of an LR-hive. Each edge signified by a
solid line within the LR-hive is assigned some fixed label. Then
these labels are sufficient to determine the labels of all the
remaining edges signified by dashed lines if either (i) any interior
edge label of the hexagon is fixed, or (ii) any boundary edge
label of the hexagon is $0$, or (iii) any two neighbouring edge
labels are equal on any of the six lines constituting the two
triangles bounding the hexagon.

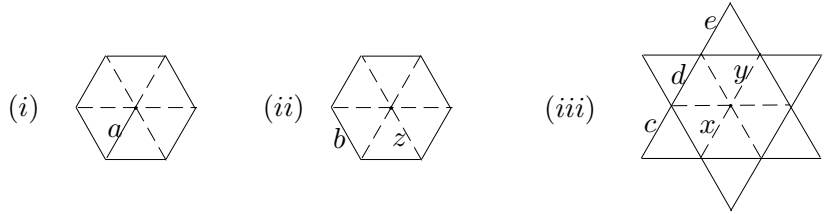
\vspace{-5mm}
\begin{center}
\input{lem4-3.TpX}
\end{center}
\end{lem}
\vspace{-5mm}

\pf The three diagrams exemplify the three possibilities referred to
in the lemma. In the first of these we can apply case (i) of
Lemma~\ref{lem-edge-labels} to each of the six elementary triangles
constituting the hexagon. These may be taken in turn, say
anticlockwise beginning with one involving the fixed edge label,
signified by $a$ in the illustrative example. In case (ii) the hive
condition~\ref{eq-hive-condition-2} gives $b\geq z$, so that for
$b=0$ we have $z=0$, since all edge labels of a LR-hive are
non-negative. Having fixed one interior edge label of the hexagon,
the remainder follow as in case (i). In case (iii) the application
of the hive conditions~\eqref{eq-hive-condition-2} gives $c\geq
x\geq d\geq y\geq e$, so that if $c=d$ we have $x=d$, and if $d=e$
we have $y=d$. In either case we have fixed one interior edge label
of the hexagon and the remainder are then fixed as in case (i). This
completes the proof, since all other examples of these three cases
can be treated in exactly the same way. \qed

Although this lemma does not exhaust the list of conditions
that fix all interior edge labels of a hexagonal subhive
of an LR-hive, avoiding these conditions turns out to be a
crucial first step in the construction of examples for which
the interior edge labels of a hexagon are not fixed. The
existence of such a situation will then be shown to characterise
those Schur function products and skew Schur functions
that are not multiplicity-free.

\section{Multiplicity-free products}
\label{sec-Stembridge}

In order to provide a hive-based proof of Stembridge's
Theorem~\ref{thm-stem} we first prove:

\begin{lem}
All the Schur function products $s_\mu s_\nu$ listed under cases
{\bf P0}--{\bf P4} of Theorem~\ref{thm-stem} are multiplicity-free.
\label{lem-Stembridge-list}
\end{lem}

\pf In order that all terms in
the product $s_\mu s_\nu$ are accounted for, we choose
$n=\ell(\mu)+\ell(\nu)$.
It then suffices to show that for
any fixed $\lambda$ there exists at most one LR-hive with boundary edge
labels specified by the parts of $\lambda$, $\mu$ and $\nu$.
This is accomplished by first parametrising the pair $\mu$ and $\nu$,
and then showing that for each fixed, but unknown
$\lambda$, the hive conditions~\eqref{eq-triangle-condition}-\eqref{eq-hive-condition-2}
serve to fix all the interior edge labels. Without the necessity of testing all
possible hive conditions, this implies that for each $\lambda$ there exists
at most one LR-hive with the required boundary edge labels, and hence that
$s_\mu s_\nu$ is multiplicity-free.

We consider the five cases in turn.

{\bf P0}.\ This case is trivial since $s_\mu\,s_0=s_\mu$ and $s_0\,s_\nu=s_\nu$
for all $\mu$ and $\nu$, respectively.

{\bf P1}.\ Thanks to the symmetry
properties~\eqref{eq-symmetry-LRcoeff}, we need only consider the
case for which $\mu=(a)$ with $a>0$ and $\nu$ fixed but arbitrary.
Then for any $\lambda$, the corresponding LR-hive takes the form in
$M1$. Applying Lemma~\ref{lem-equal-edge-labels} to the triangle
$BCD$ fixes all its edge labels, including those on $BD$. Case (iv)
of Lemma~\ref{lem-edge-labels} then serves to fix all the edge
labels of $ABDE$. Thus all edge labels of the complete LR-hive are
fixed, so that each product of type {\bf P1} is multiplicity-free.

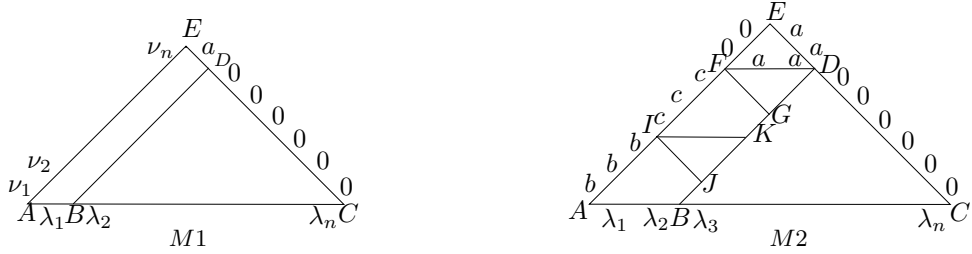
\vspace{-5mm}
\begin{center}
\input{mpf-stem-12.TpX}
\end{center}
\vspace{-5mm}

{\bf P2}.\ Thanks once again to the symmetry
properties~\eqref{eq-symmetry-LRcoeff}, we need only consider the
case for which $\mu=(a^2)$ with $a>0$ and $\nu=(b^pc^q)$ with
$b>c>0$, $p,q>0$ and $n=p+q+2$. If $p,q>2$ then, for any $\lambda$
with $\ell(\lambda)\leq n$, the corresponding LR-hives take the form
in $M2$. Lemma~\ref{lem-equal-edge-labels} implies that all the edge
labels of $BCD$, $ABJI$ and $DEF$ are fixed. Since those on $DF$
must both be $a$, Lemma~\ref{lem-equal-edge-labels} implies that all
the edge labels of $DFG$ are also fixed. This is then sufficient to
determine all edge labels of $FIKG$, thanks again to
Lemma~\ref{lem-equal-edge-labels}. This only leaves the triangle
$IJK$ of side length $2$ to be considered. Its boundary edges labels
are all known, so that, by virtue of case (ii) of
Lemma~\ref{lem-edge-labels}, its interior edge labels are also
fixed. This fixes all edge labels in the complete LR-hive, and the
corresponding product is multiplicity-free.

If either $q=2$ or $p=2$, then the previous diagram must be modified
as shown below. The argument then proceeds exactly as before with
the trapeziums $ABJI$ and $FIKG$ replaced by the triangles $AJI$ and
$FIK$, as appropriate.

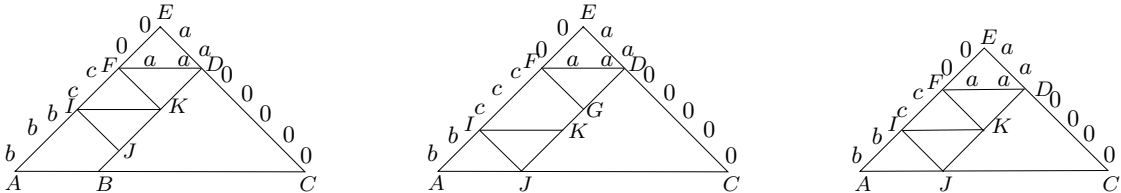
\vspace{-5mm}
\begin{center}
\input{mpf-stem-21.TpX}
\end{center}
\vspace{-5mm}

On the other hand if either $p=1$ or $q=1$, then $\nu$ is a near
rectangle, and this is a situation covered by case {\bf P3}.

{\bf P3}.\ Again thanks to the symmetry
properties~\eqref{eq-symmetry-LRcoeff} it is sufficient to consider
the two cases (i) $\mu=(a^pb)$ and $\nu=(c^q)$ and (ii) $\mu=(ab^p)$
and $\nu=(c^q)$ with $a>b>0$, $c,p,q>0$ and $n=p+q+1$. In each
subcase (i) and (ii) there are three possibilities, depending on the
relative size of $p$ and $q$. For any $\lambda$, with
$\ell(\lambda)\leq n$, the LR-hives in the subcase (i) may take one
or other of the following forms:

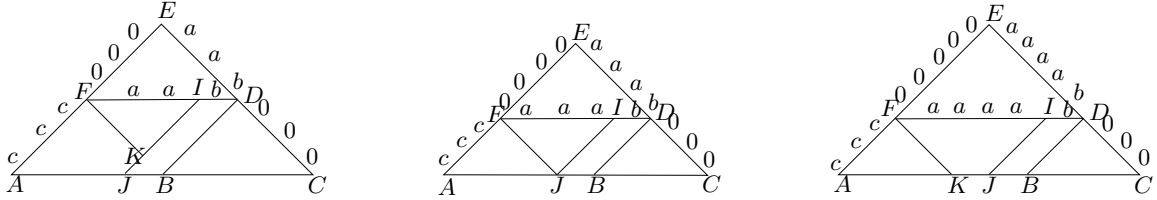
\vspace{-5mm}
\begin{center}
\input{mpf-stem-4i.TpX}
\end{center}
\vspace{-5mm}

Considering the first diagram, successive applications of
Lemma~\ref{lem-equal-edge-labels} fix all the edge labels of the
triangles $DEF$ and $BCD$, as well as those of the trapezium $AJKF$.
Since the edge labels of $FI$ are all $a$, all the edge labels of
$FIK$ are then fixed by virtue of Lemma~\ref{lem-equal-edge-labels}.
This means that the boundary edge labels of the thin strip $BDIJ$
are all known, so that thanks to case (iii) of
Lemma~\ref{lem-edge-labels}, all the remaining interior edge labels
are also fixed. This completes the edge labelling of the complete
hive. A similar argument applies to the other two diagrams.

Similarly, for the subcase (ii), we have the following types of
LR-hive, and the argument goes through precisely as before.

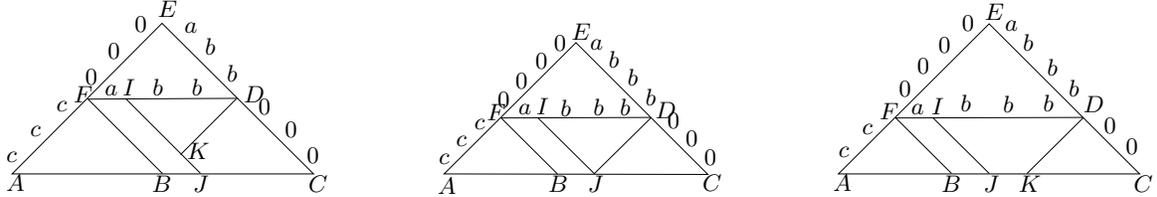
\vspace{-5mm}
\begin{center}
\input{mpf-stem-4ii.TpX}
\end{center}
\vspace{-5mm}

It follows that the case {\bf P3} is also multiplicity-free.

{\bf P4}.\ For this case, let $\mu=(a^p)$ and $\nu=(b^q)$ with
$a,b,p,q>0$ and $n=p+q$. There are three subcases corresponding to
$p<q$, $p=q$ and $p>q$. For each of these, for any $\lambda$ with
$\ell(\lambda)\leq n$, the corresponding LR-hives take the form:

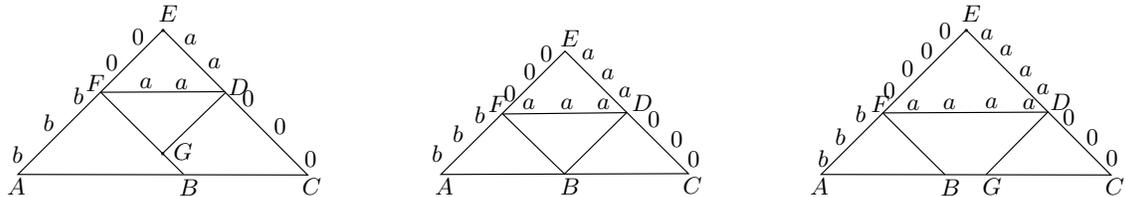
\vspace{-5mm}
\begin{center}
\input{mpf-stem-3.TpX}
\end{center}
\vspace{-5mm}

In the first of these LR-hives, thanks to
Lemma~\ref{lem-equal-edge-labels} all the edge labels of the
triangles $DEF$ and $ABF$, as well as those of the trapezium $BCDG$,
are fixed. Since the edge labels on $FD$ are all $a$,
Lemma~\ref{lem-equal-edge-labels} fixes all the edge labels of
$DFG$, thereby completing the edge labelling of the complete hive.
Thus at most one LR-hive exists. A similar argument applies to the
other two diagrams, and this case {\bf P4} is also
multiplicity-free.

This completes the proof of Lemma~\ref{lem-Stembridge-list}. \qed

\section{Completeness of the list in Stembridge's theorem}
\label{sec-Stembridge-completeness}

To complete the proof of Stembridge's Theorem~\ref{thm-stem} it is
necessary to show that all cases other than those of {\bf P0}--{\bf
P4} are not multiplicity-free. This can be done by following
Stembridge's argument based on the use of just the second part
of~\eqref{eq-add-col-row}. In the context of the hive model, we do
this by first considering three further cases, for which we shall
show that there exists at least one partition $\lambda$ such that
$c^{\lambda}_{\mu\nu}>1$.

\begin{lem}
The product $s_{\mu}s_{\nu}$ is not multiplicity-free in each of
the following cases:
\begin{description}
\item[Q1]~~~~$\mu=(ab)$ and $\nu=(cd)$ with $a>b>0$, $c>d>0$;

\item[Q2]~~~~$\mu=(abc)$ and $\nu=(d^2)$ with $a>b>c>0$ and $d>1$;

\item[Q3]~~~~$\mu=(a^2b^2)$ and $\nu=(c^3)$ with $a>b+1$, $b>1$ and $c>2$.
\end{description}
\label{lem-stem}
\end{lem}

\pf For $\lambda=(xyz)$, $(wxyz)$ and $(uvxyzw)$ the corresponding
hives take the form shown below, with $p=a+d-w$ in $H2$,
and $p=a+c-v$, $q=a+c-u$ and $r=a+b+c-u$ in $H3$:

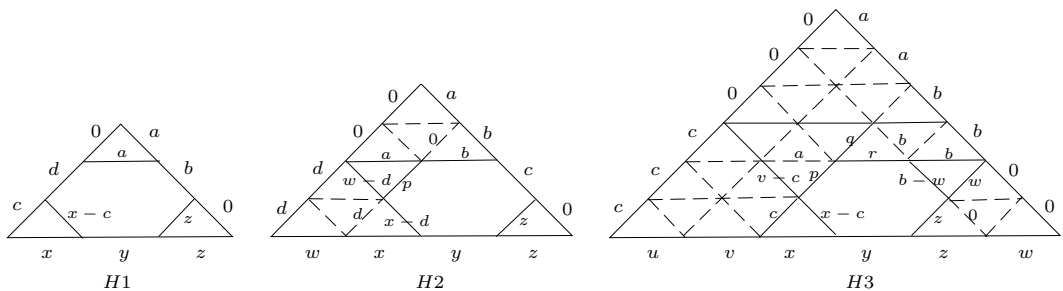
\vspace{-5mm}
\begin{center}
\input{lem5-1a.TpX}
\end{center}
\vspace{-5mm}

In each case the solid lines divide the hive into portions for which
the edge labels are determined, including all the dashed line
interior edges. In each case it will be observed that we are left
with a hexagon on which the six boundary edge labels are necessarily
fixed from a knowledge of $\lambda$, $\mu$ and $\nu$.

Now, for each given $\mu$ and $\nu$ we will identify one particular partition
$\lambda$ for which there exists exactly two distinct labellings of
the interior edges of the hexagon that satisfy all the hive
conditions. It will then follow that $c^\lambda_{\mu\nu}=2$, so
that $s_{\mu}s_{\nu}$ is not multiplicity-free.

We consider each case in turn.

\noindent{\bf Q1}.~~
If we take $\lambda=(a+c-1,b+d,1)$, then with the stated conditions,
there are exactly two LR-hives $H1$ corresponding to this given $\lambda$,
as the following figures show:

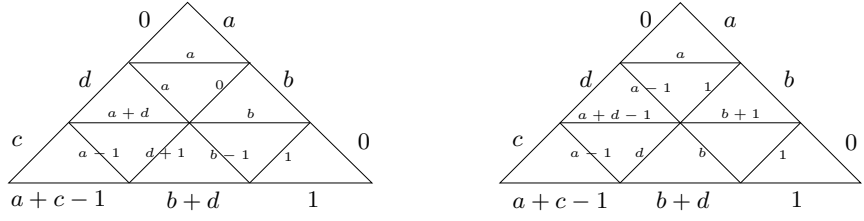
\vspace{-5mm}
\begin{center}
\input{ab-cd.TpX}
\end{center}
\vspace{-5mm}

Thus $s_{ab}s_{cd}$ is not multiplicity-free.

\noindent{\bf Q2}.~~
For $\lambda=(a+d-1,b+d-1,c+1,1)$, with the stated conditions,
we can complete the
labelling of the interior edges of two LR $3$-subhives of $H2$, as
illustrated in the following figures:

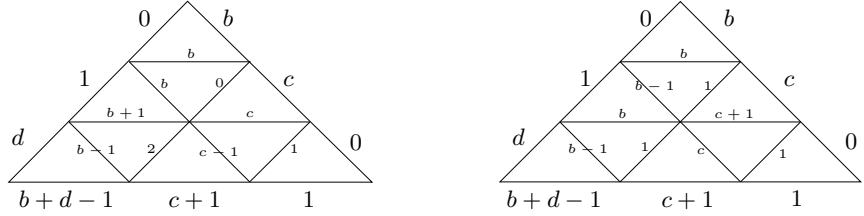
\vspace{-5mm}
\begin{center}
\input{abc-dd.TpX}
\end{center}
\vspace{-5mm}

The edge labellings of these pairs of LR $3$-hives serve to complete
the interior edge labelling of the corresponding pairs of LR
$4$-hives in which they are embedded.
The existence of two LR-hives
corresponding to the given $\lambda$ shows that $s_{abc}s_{d^2}$ is
not multiplicity-free.

\noindent{\bf Q3}.~~
For $\lambda=(a+c-1,a+c-2,b+c-1,b+1,2,1)$, with the stated conditions,
we can complete the labelling of the interior edges of two LR $3$-subhives
of $H3$ as illustrated in the following figures:

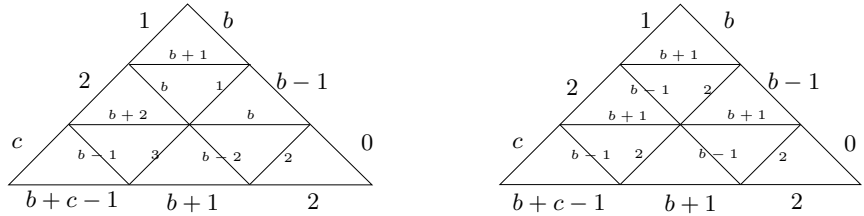
\vspace{-5mm}
\begin{center}
\input{aabb-ccc.TpX}
\end{center}
\vspace{-5mm}

The edge labellings of these pairs of LR $3$-hives serve to complete
the interior edge labelling of the corresponding pairs of LR
$6$-hives in which they are embedded.
The existence of two LR-hives
corresponding to the given $\lambda$ shows that $s_{a^2b^2}s_{c^3}$
is not multiplicity-free.

This completes the proof of Lemma \ref{lem-stem}.\qed

\noindent{\bf Note}\ \
It should be pointed out that the conditions
on $a, b, c, d$ as stated in Lemma \ref{lem-stem} arise naturally.
In $H1$, to avoid being multiplicity-free,
Lemma~\ref{lem-hexagon} implies that $a>b>0$ and $c>d>0$, as stated
for case {\bf Q1}. In $H2$, Lemma~\ref{lem-hexagon} implies that we
require $a>b$, $b>c>0$ and $d>p>0$ so that $a>b>c>0$ and $d>1$ as
stated for case {\bf Q2}. In $H3$, to avoid being
multiplicity-free, Lemma~\ref{lem-hexagon} gives $a>r>b$, $b>b-w>0$
and $c>p>q$. Hence $a>b+1$, $b>1$ and $w>0$. Then, in order to avoid
fixing an interior edge of the hexagon by means of the hive
condition~\eqref{eq-hive-condition-2} we must also have $p>w$, so
that $c>p>w>0$. Hence $c>2$, as required to complete the conditions
listed for case {\bf Q3}.

\noindent{\bf Proof of Theorem \eqref{thm-stem}}. Suppose $\mu$ and
$\nu$ have $s$ and $t$ distinct non-zero parts, respectively, with
$s,t\geq 0$. For $s>0$ and $t>0$ we 
let $\mu=(a_1^{p_1},a_2^{p_2},\ldots,a_s^{p_s})$ and
$\nu=(b_1^{q_1},b_2^{q_2},\ldots,b_t^{q_t})$, with
$a_1>a_2>\cdots>a_s>0$ and $b_1>b_2>\cdots>b_t>0$, where $p_i, q_j>0$
for $i=1,2,\ldots,s$ and $j=1,2,\ldots,t$.
First we recall that the results of Section~\ref{sec-Stembridge}
imply that $s_{\mu}s_{\nu}$ is multiplicity-free in each
of the cases {\bf P0}--{\bf P4}. Then
we consider all possible values of $s$ and $t$ in turn.

In the following, we select the parts of $\sigma$ and $\tau$
from those of $\mu$ and $\nu$, respectively, so that
$\mu=\sigma\cup\zeta$ and $\nu=\tau\cup\xi$ for some
$\zeta$ and $\xi$. If the choice is made in such a way that
$s_{\sigma}s_{\tau}$ is not multiplicity-free,
then there exists at least one $\rho$ such
that $c_{\sigma\tau}^{\rho}\geq 2$.
Now let $\lambda=\rho\cup\eta$ where $\eta$ is formed
by pairing up the parts of $\zeta$ and $\xi$ in any convenient way
so that each $\eta_k=\zeta_i+\xi_j$ for some $i$ and $j$. It then
follows from the repeated application of the second part
of~\eqref{eq-add-col-row} with $a=\eta_k$, $b=\zeta_i$ and $c=\xi_j$
that $c^\lambda_{\mu\nu}\geq c^\rho_{\sigma\tau}\geq2$. Thus
$s_{\mu}s_{\nu}$ is not multiplicity-free.

The cases $s=0$ or $t=0$ are covered by case {\bf P0}, and 
are multiplicity-free. 
If $s\geq2$ and $t\geq2$ we select
$\{a,b\}\subseteq\{a_1,a_2,\ldots,a_s\}$ and
$\{c,d\}\subseteq\{b_1,b_2,\ldots,b_t\}$ in such a way that
$\sigma=(ab)$ and $\tau=(cd)$ are a pair of partitions of the type
covered by case {\bf Q1}.

If $s\geq3$ the case $t\geq2$ has already been dealt with. We can
therefore take $t=1$ so that $\nu=(b^q)$. If $q=1$ or $b=1$ then
$s_{\mu}s_{\nu}$ is multiplicity-free since the situation is covered
by case {\bf P1}, which were shown to be multiplicity-free before.
For $q>1$ and $b>1$, it is possible to select from the distinct
parts of $\mu$ and $\nu$ those parts that constitute partitions
$\sigma=(abc)$ and $\tau=(d^2)$ appropriate to the
non-multiplicity-free case {\bf Q2}.

If $s=2$ and $t=1$, suppose $\mu=(a^pb^q)$ and $\nu=(c^r)$. If $c=1$
or $r=1$, then the situation is covered by case {\bf P1}; if $c=2$
or $r=2$, then the situation is covered by case {\bf P2}; if $a=b+1$
or $b=1$ or $p=1$ or $q=1$, then the situation is covered by case
{\bf P3}. For all these cases $s_{\mu}s_{\nu}$ is multiplicity-free.
Thus we consider the case where $a>b+1$, $b>1$, $p,q>1$ and $c,r>2$.
In this case we can always select $\sigma=(a^2b^2)$ and
$\tau=(c^3)$, and this is covered by the non-multiplicity-free case
{\bf Q3}.

Any case with $s=1$ and $t\geq 2$ is related by the first symmetry
condition of ~\eqref{eq-symmetry-LRcoeff} to a case with $t=1$ and
$s\geq 2$ that has already been dealt with, so only the case $s=1$
and $t=1$ is left. This case appears in the list under {\bf P4}, and
is multiplicity-free.

This completes the proof of Theorem~\ref{thm-stem}. \qed

\section{Multiplicity-free skew Schur functions}
\label{sec-main-result}

In this section we prove:

\begin{lem}
All the basic skew Schur functions $s_{\lambda/\mu}$ listed under
cases {\bf R0}--{\bf R4} of Theorem~\ref{thm-main} are multiplicity-free.
\label{lem-main}
\end{lem}

\pf For each case the strategy is to parametrise the pair
$\lambda$ and $\mu$ and then, for each fixed but unknown $\nu$, to
use the hive conditions
\eqref{eq-triangle-condition}--\eqref{eq-hive-condition-2} to show
that all the interior edge labels are fixed.
Without the necessity of testing all
possible hive conditions, this implies that for each $\nu$ there
exists at most one LR-hive with the required boundary edge labels,
and hence that $s_{\lambda/\mu}$ is multiplicity-free.

Since $s_{\lambda/\mu}$ is basic, the required LR-hives are integer
$n$-hives, with $n=\ell(\lambda)$ and all edge
labels positive along the boundary specified by $\lambda$ and at
least one edge label $0$ along the boundary specified by $\mu$.
{\it We consider the four %%%eight
cases in turn.

{\bf R0}.~~There are two subcases, {\bf S0}: $\mu=0$, and ${\bf
S0^\pi}$: $\lambda=(m^n)$ for some positive integer $m$, as
illustrated by: 
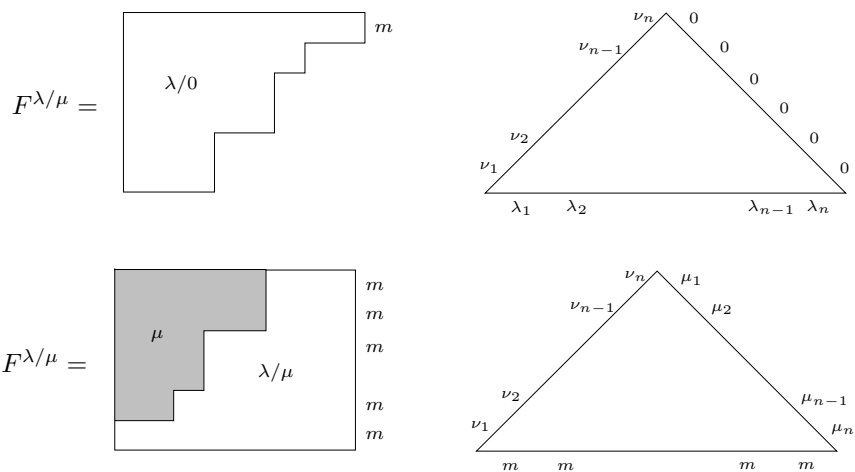
\vspace{-3mm}
\begin{center}
\input{4-1-1.TpX}

\input{4-1-2.TpX}
\end{center}
\vspace{-5mm} Lemma~\ref{lem-equal-edge-labels} implies immediately
that in each case there exists a single LR-hive. In the case of {\bf
S0} the equal edge labels $0$ suffice to show that $\nu=\lambda$,
while in the case of ${\bf S0^\pi}$ the equal edge labels $m$ fix
the parts of $\nu$ to be $\nu_k=m-\mu_{n-k+1}$ for $k=1,2,\ldots,n$.

{\bf R1}.~~There are two major subcases which we designate by {\bf
S1} and ${\bf S1^\pi}$ in which $\mu$ and $\lambda^*$, respectively,
are rectangles of $m^n$-shortness $1$. Each has four subcases, as
illustrated by:

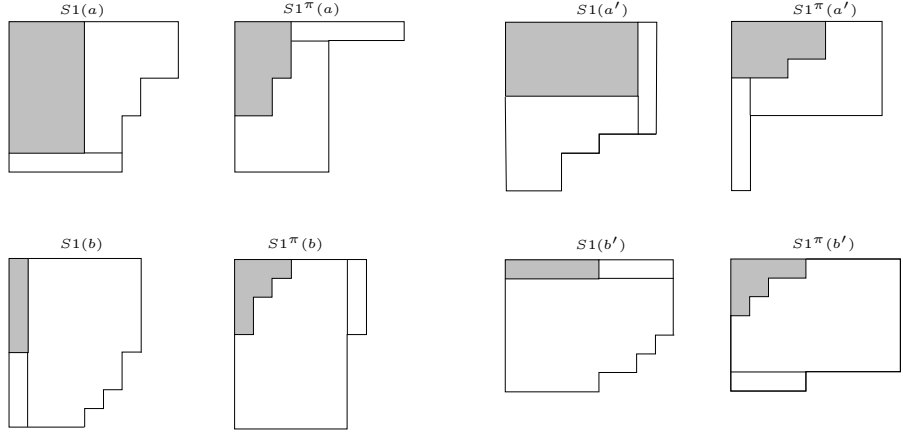
\begin{center}
\input{4-12-5new.TpX}
\end{center}
\vspace{-5mm}

\noindent These skew Young diagrams for $s_{\lambda/\mu}$ have been
arranged so that those of type ${\bf S1^{\pi}}$ are just the
$\pi$-rotations of {\it their left-hand neighbour} of type {\bf S1}.
Moreover, the right-hand block of four are just the conjugates of
the left-hand block. Thanks to the rotation
symmetry~\eqref{eq-skew-rotate} and the conjugate
symmetry~\eqref{eq-symmetry-LRcoeff}, it is therefore only necessary
to consider two cases, which we choose to be ${\bf S1^{\pi}(a)}$ and
${\bf S1^{\pi}(b')}$.}

${\bf S1^{\pi}(a)}$. Suppose $\lambda=(ab^{n-1})$ and $\mu$ is
arbitrary, then the corresponding Young diagram and LR-hives take
the form:

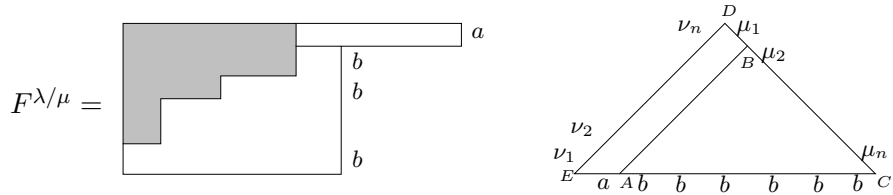
\vspace{-5mm}
 \begin{center}
 \input{4-8.TpX}
 \end{center}
 \vspace{-5mm}

 \noindent For given $\mu$ and $b$,
Lemma~\ref{lem-equal-edge-labels} implies that all the edge labels
of $ABC$ are fixed, including those on $AB$. It then follows from
case (iv) of Lemma~\ref{lem-edge-labels} that, for any given $\nu$,
all the edge labels of $ABDE$ are also fixed. Thus all the hive edge
labels are fixed, and $s_{\lambda/\mu}$ must be multiplicity-free,
as required.

${\bf S1^{\pi}(b')}$. Suppose $\lambda=(a^{n-1}b)$ and $\mu$ is
arbitrary. Then this case is exemplified by:

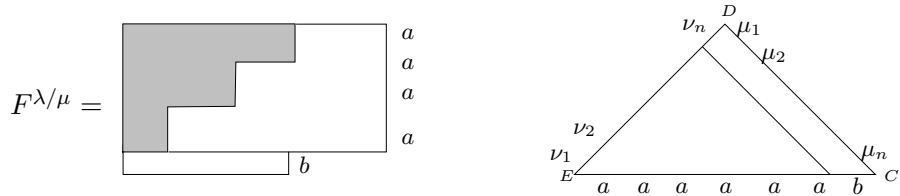
\vspace{-5mm}
\begin{center}
\input{4-9.TpX}
\end{center}
\vspace{-5mm}

\noindent That $s_{\lambda/\mu}$ is multiplicity-free then follows
from an argument entirely analogous to that used for ${\bf
S1^{\pi}(a)}$.

{\it This completes the argument that each $s_{\lambda/\mu}$ of
{\bf R1} is multiplicity-free.
}

%%%%%%%%%%%%%%%%%%%%%
{\it

{\bf R2}.~~The two major subcases, $\mu$ a rectangle of
$m^n$-shortness $2$ and $\lambda^*$ a fat hook, and {\it vice
versa}, we designate by {\bf S2} and ${\bf S2^{\pi}}$, respectively.
They each possess four subcases, as illustrated by:

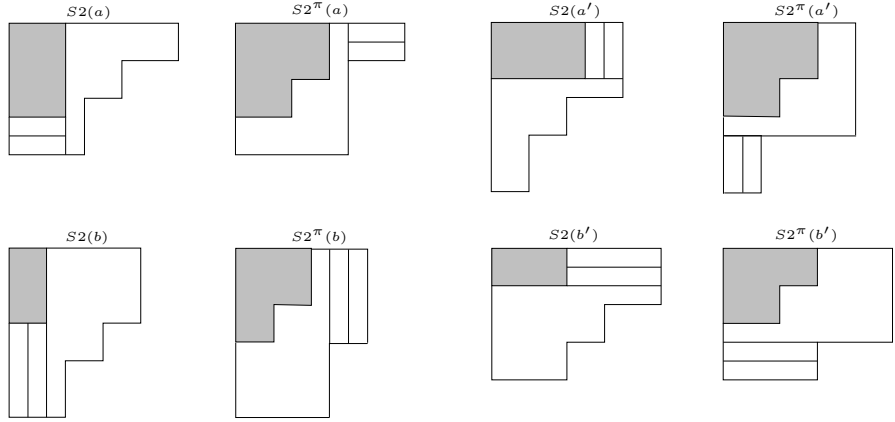
\begin{center}
\input{6-8new.TpX}
\end{center}

\noindent Once again, these skew Young diagrams have been arranged
so that those of type ${\bf S2^{\pi}}$ are just the $\pi$-rotations
of those of type {\bf S2}, and the right-hand block of four is just
the conjugate of the left-hand block of four. Thanks to the rotation
symmetry~\eqref{eq-skew-rotate} and the conjugate
symmetry~\eqref{eq-symmetry-LRcoeff}, it is therefore only necessary
to consider two cases, which we choose to be {\bf S2(a)} and ${\bf
S2(b')}$. }

{\bf S2(a)}.~~In this case $\lambda=(a^r b^s c^t)$ and
$\mu=(d^p0^q)$ with $r+s+t=p+q=n$ and $q=2$, as illustrated in the
following figure:

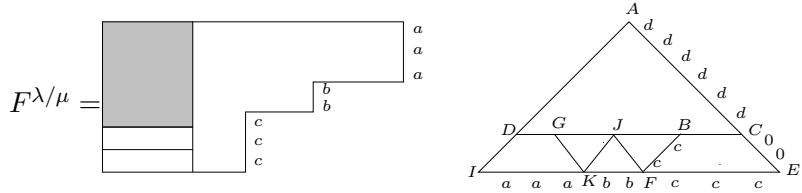
\vspace{-5mm}
\begin{center}
\input{6-7c.TpX}
\end{center}
\vspace{-50mm}

First, all the edge labels in subhives $ACD$, $IKGD$ and $CEF$ can
be determined uniquely by Lemma \ref{lem-equal-edge-labels} with all
the edge labels on $CF$ being $c$, by
Lemma~\ref{lem-equal-edge-labels} again, the labels in $BCF$ can be
determined uniquely and then all the edge labels in $KFBJ$ can be
determined uniquely. Since the edge labels on the boundary of $GKJ$
are known, and this triangle has side length two, then all its
interior edge labels are fixed by those on the boundary by using
case (ii) of Lemma \ref{lem-edge-labels}. Hence all edge labels are
fixed and this case is also multiplicity-free.

${\bf S2(b')}$.~~In this case $\lambda=(a^r b^s c^t)$ and
$\mu=(d^p0^q)$ with $r+s+t=p+q=n$ and $p=2$, as illustrated in the
following figure:

\vspace{-5mm}
\begin{center}
\input{4-20.TpX}
\end{center}
\vspace{-45mm}

 \noindent By Lemma~\ref{lem-equal-edge-labels}, all
the edge labels in $ABC$, $IKH$ and $CJKE$ can be determined
uniquely, and the edge labels on $JF$ and $FC$ are all $b$ and $c$,
respectively. It follows, that all the edge labels in the subhives
$BCF$ and $HJFG$ can also be determined by
Lemmas~\ref{lem-equal-edge-labels} again. This only leaves the
interior labels of $BGF$ undetermined, which can be determined
immediately by case (ii) of Lemma \ref{lem-edge-labels}. Then all
the edge labels in the complete hive are determined, and once again
this case is multiplicity-free.

The symmetry conditions~\eqref{eq-skew-rotate}
and~\eqref{eq-symmetry-LRcoeff} then establish the fact that all
case {\bf R2} examples are multiplicity-free.

%%%%%%%%%%%%%%%%%%%%%
{\it

{\bf R3}.~~The two major subcases, $\mu$ a rectangle and $\lambda^*$
a fat hook of $m^n$-shortness $1$, and {\it vice versa}, we
designate by {\bf S3} and ${\bf S3^{\pi}}$, respectively. They each
possess six subcases, as illustrated by:

\begin{center}
\input{5anew.TpX}
\end{center}

\noindent Once again, these skew Young diagrams have been arranged
so that those of type ${\bf S3^{\pi}}$ are just the $\pi$-rotations
of those of type {\bf S3}. This time the right-hand block of six is
just the conjugate of the left-hand block of six. Thanks to the
rotation symmetry~\eqref{eq-skew-rotate} and the conjugate
symmetry~\eqref{eq-symmetry-LRcoeff}, it is therefore only necessary
to consider three cases, which we choose to be {\bf S3(a)}, {\bf
S3(b)} and {\bf S3(c)}. }

{\bf S3(a)}. Suppose $\lambda=(ab^sc^t)$, $\mu=(d^p0^q)$ with
$s,t,p,q>0$ and $1+s+t=p+q=n$. The case $p=1$ and $p=s+t$ have been
covered in ${\bf S1(b')}$ and {\bf S1(a)}. This leaves three cases
to discuss: $1<p<1+s$, $p=1+s$, and $1+s<p<s+t$, as the following
figure shows:

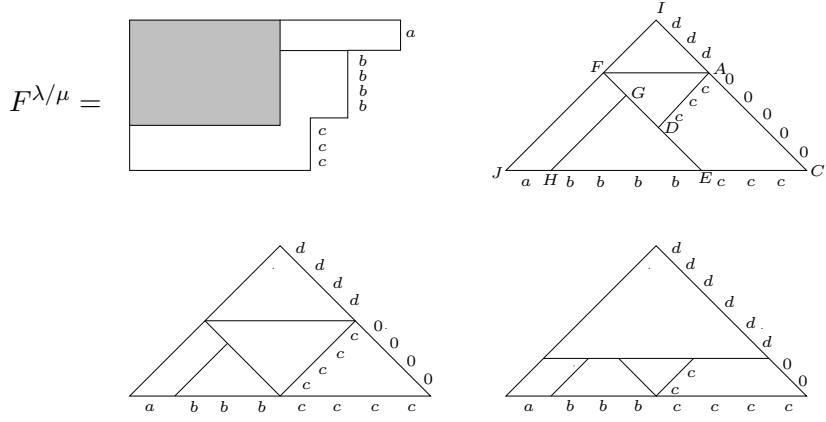
\vspace{-5mm}
\begin{center}
\input{4-15.TpX}
\end{center}
\vspace{-5mm}

We only consider the first case $1<p<1+s$. The argument for the others is similar.

First, by Lemma~\ref{lem-equal-edge-labels}, the edge labels in
$AIF$ and $ADEC$ can be determined, and the edge labels on $AD$ are
all equal to $c$. Again Lemma~\ref{lem-equal-edge-labels} implies
that all edge labels in $ADF$ are fixed. Thus all edge labels on
$GE$ are determined. The use, yet again, of
Lemma~\ref{lem-equal-edge-labels} suffices to fix all edge labels in
$HEG$. Finally, the edge labels of the region $FJHG$ are fixed by
virtue of case (iv) of Lemma~\ref{lem-edge-labels}. Thus we have
determined all the edge labels of this hive. Hence $s_{\lambda/\mu}$
is multiplicity-free.

{\bf S3(b)}. Suppose $\lambda=(a^rbc^t)$, $\mu=(d^p0^q)$ with
$r,t,p,q>0$ and $r+1+t=p+q=n$. Since $p=1$ and $p=r+t$ have been
covered in ${\bf S1(b')}$ and {\bf S1(a)}, there are four subcases
$1<p<r$, $p=r$, $p=r+1$ and $r+1<p<r+t$, as illustrated below:

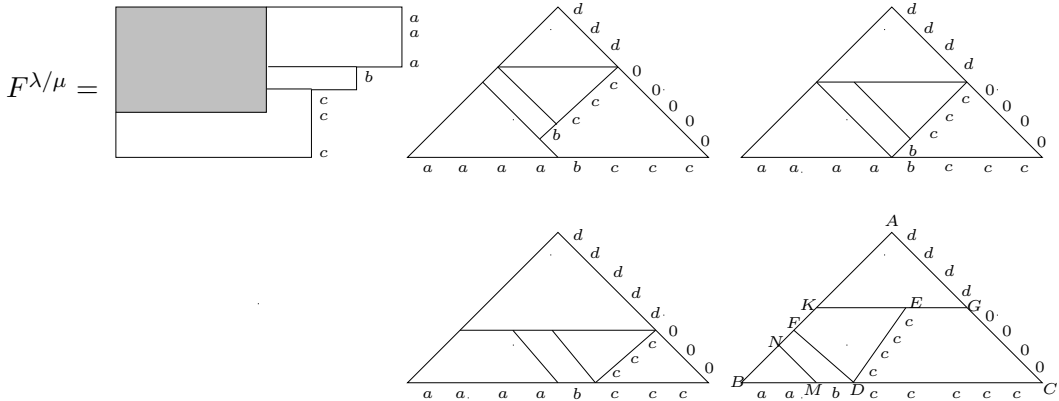
\vspace{-5mm}
\begin{center}
\input{4-17.TpX}
\end{center}
\vspace{-5mm}

 \noindent By way of example, we consider the fourth
subcase. The others may be dealt with similarly. By
Lemma~\ref{lem-equal-edge-labels}, all the edge labels in subhives
$BMN$, $AKG$ and $DCGE$ are fixed, and the edge labels on $DE$ are
all $c$. Then by Lemma~\ref{lem-equal-edge-labels} once again, the
edge labels in $DEKF$ can be determined. Finally, thanks to
Lemma~\ref{lem-edge-labels} the edge labels of the thin strip $MDFN$
are completely determined. We can therefore conclude that
$s_{\lambda/\mu}$ is again multiplicity-free.

{\bf S3(c)}. Suppose $\lambda=(a^rb^sc)$ and $\mu=(d^p0^q)$ with
$r,s,p,q>0$ and $r+s+1=p+q=n$. Since $p=1$ and $p=r+s$ have been
covered in ${\bf S1(b')}$ and {\bf S1(a)} respectively, there are
three cases to consider: $1<p<r$, $p=r$ and $r<p<r+s$. We choose to
illustrate just the case $r<p<r+s$:

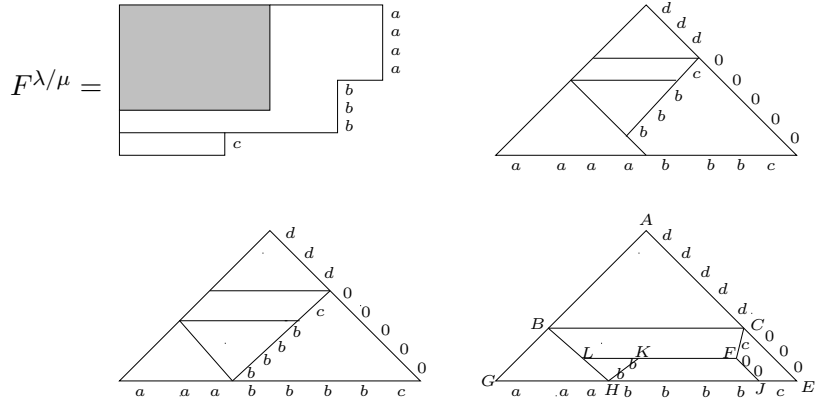
\vspace{-5mm}
\begin{center}
\input{4-16.TpX}
\end{center}
\vspace{-5mm}

 \noindent By Lemma~\ref{lem-equal-edge-labels}, the
edge labels in $ABC$ and $BGH$ are fixed by the hive boundary edge
labels and the edge labels $0$ along $CE$ force all the edge labels
on $FJ$ to also be $0$. Thanks to Lemma~\ref{lem-equal-edge-labels}
it follows that the edge labels in $HJFK$ are fixed, with those on
$HK$ all equal to $b$. Then Lemma~\ref{lem-equal-edge-labels} fixes
all the edge labels in $LHK$. This leaves the thin strips $BLFC$ and
$CFJE$ each of which may be dealt with through the use of
Lemma~\ref{lem-edge-labels}. This fixes all the edge labels in the
complete hive. Once again $s_{\lambda/\mu}$ is multiplicity-free.
Similar arguments cover all the other subcases.

Hence, by the symmetry conditions~\eqref{eq-skew-rotate} and
~\eqref{eq-symmetry-LRcoeff}, it follows that all the skew Schur
functions of case {\bf R3} are also multiplicity-free.

%%%%%%%%%%%%%%%%%%%%%%%%
{\it

{\bf R4}.~~Here both $\mu$ and $\lambda^*$ are rectangles. The
latter implies that $\lambda$ itself is a fat hook. In this
situation, which we designate by {\bf S4}, } We can set
$\lambda=(a^rb^s)$, $\mu=(c^p0^q)$ with $r,s,p,q>0$ and $r+s=p+q=n$.
There are three subcases to cover: $p<r$, $p=r$ and $p>r$. These are
illustrated by:

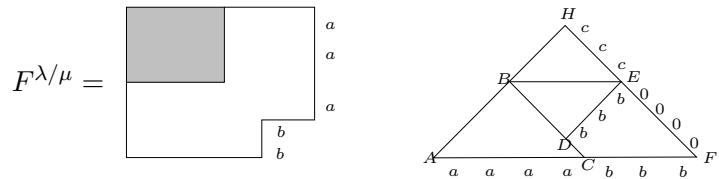
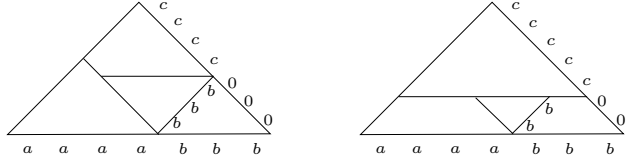
\vspace{-5mm}
\begin{center}
\input{4-13.TpX}

\input{4-14.TpX}
\end{center}
\vspace{-5mm}

\noindent  In the first of these, by
Lemma~\ref{lem-equal-edge-labels}, the edge labels in $ABC$, $BEH$
and $CFED$ can be determined and the edge labels on $DE$ are all
equal to $b$. Since, in addition, the edge labels on $BD$ are known,
all the edge labels in $BDE$ are fixed by Lemma
\ref{lem-equal-edge-labels}. Thus we have determined all the edge
labels of this hive. A similar argument applies in the other two
subcases. Thus in all three subcases $s_{\lambda/\mu}$ is
multiplicity-free.

%%%%%%%%%%%%%%%%%%%%%%

This completes the proof of Lemma~\ref{lem-main} that all the cases
listed in Theorem~\ref{thm-main} are multiplicity-free, as claimed.
\qed

\section{Completeness of the list in the main theorem}
\label{sec-main-completeness}

It remains to show that the list of multiplicity-free skew Schur
functions given in Theorem~\ref{thm-main} is exhaustive. To this end
we first consider three further cases, for which we shall show that
there exists at least one partition $\nu$ such that
$c^\lambda_{\mu\nu}>1$.

\begin{lem} The skew Schur functions $s_{\lambda/\mu}$ are
not multiplicity-free in each of the following cases:
\begin{description}
\item[T1]~~~~$\lambda=(abc)$ and $\mu=(de)$ with $a>b>c>0$, $d>e>0$, $a>d$ and $b>e$;

\item[T2]~~~~$\lambda=(abcd)$ and $\mu=(e^2)$ with $a>b>c>d>0$ and $b>e>1$;

\item[T3]~~~~$\lambda=(a^2b^2c^2)$ and $\mu=(d^3)$ with $a>b+1$, $b>c+1$, $c>1$ and $b>d>2$.
\end{description}
\label{lem-multiplicity-main}
\end{lem}

\pf\

For $\nu=(xyz)$, $(xyzw)$ and $(wxyzuv)$ the
corresponding LR-hives take the form shown below,
with $p=e-d+w$ in $K2$, and
$p=a+b-w$, $q=b+c-d-v$, $r=d-c+u$ and $s=d-c+v$ in $K3$:

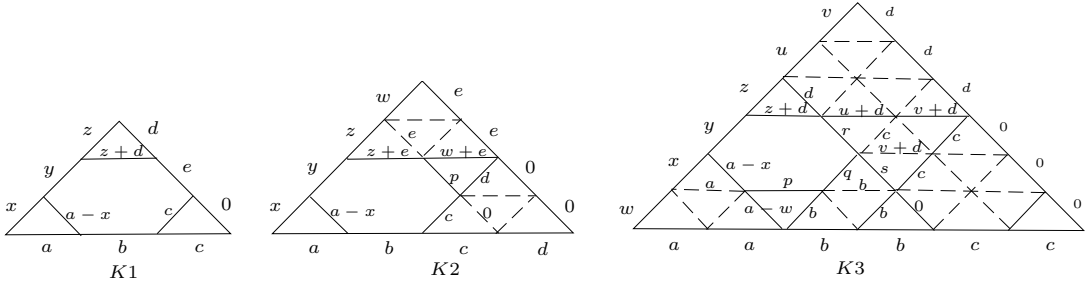
\begin{center}
\input{lem7-1a.TpX}
\end{center}
\vspace{-7mm}

The solid lines divide the hive into portions for which the edge
labels are determined, including all the dashed line interior edges.
In each case, we are left with a hexagon on which the six boundary
edge labels are necessarily fixed from a knowledge of $\lambda$,
$\mu$ and $\nu$.

A priori the skew Schur functions $s_{\lambda/\mu}$
identified here need not be basic.
They are row-basic since the stated conditions
ensure that $\lambda_i>\mu_i$ for all $i$, so that no row of
$F^{\lambda/\mu}$ is empty. However, some columns may be empty.
If so, then these may be deleted to give $F^{\hat\lambda/\hat\mu}$,
with $s_{\hat\lambda/\hat\mu}$ basic.
Since the pair $\hat\lambda$ and $\hat\mu$
belong to the same case
{\bf T1}, {\bf T2} or {\bf T3} as the original pair $\lambda$
and $\mu$, and $s_{\lambda/\mu}=s_{\hat\lambda/\hat\mu}$,
in accordance with \eqref{eq-skew-basic}, we can,
without loss of generality, confine our attention to
those cases for which $s_{\lambda/\mu}$ is itself basic.

Then, for each given pair $\lambda$ and $\mu$,
such that $s_{\lambda/\mu}$ is basic, we will identify one
particular partition $\nu$ for which there are
precisely two distinct labellings of
the interior edge labels of the hexagon that satisfy all the hive
conditions. It will then follow that $c^\lambda_{\mu\nu}=2$, so that
$s_{\lambda/\mu}$ is not multiplicity-free. We consider each case in turn.

\noindent{\bf T1}.~~Here, with $\lambda=(abc)$ and $\mu=(de)$
and $s_{\lambda/\mu}$ basic,
{\it there are just two overlapping subcases to
consider, one with $c+1\geq d$ and the other with $d\geq c+1$.
In each subcase, we offer an appropriate partition $\nu=(xyz)$ for
which there exist two LR-hives:}
$$
\begin{array}{lclcl}
{\bf T1}{\rm(i)}&&a>b\geq c+1\geq d\geq e+1>1&&\nu=(a-1,b-e,c-d+1)\,;\cr
{\bf T1}{\rm(ii)}&&a>b\geq d\geq c+1\geq e+1>1&&\nu=(a-1,b+c-d-e+1,0)\,.\cr
\end{array}
$$
For these $\nu$ the corresponding pairs of LR-hives are given explicitly by:

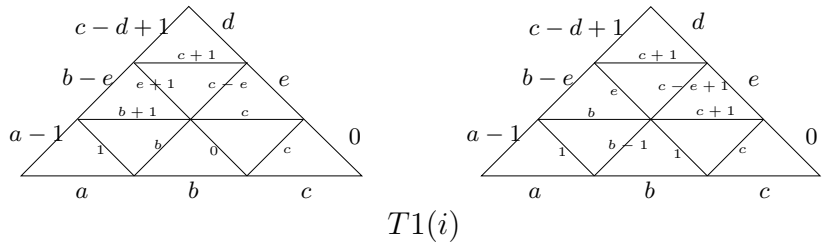
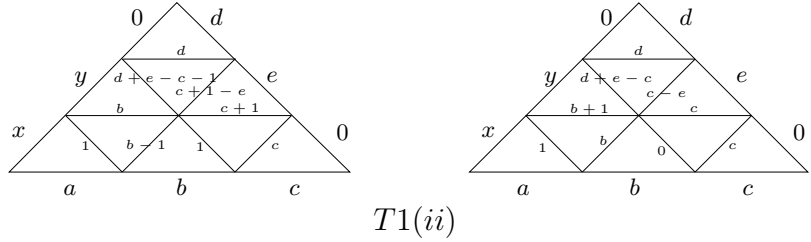
\begin{center}
\input{abc-de05.TpX}
\end{center}

\begin{center}
\input{abc-de04.TpX}
\end{center}

\vspace{-5mm}

Thus there are precisely two LR-hives corresponding to the given
$\nu$, so that $s_{abc/de}$ is not multiplicity-free.
It might be
noted that the case $a=3$, $b=2$, $c=1$, $d=2$ and $e=1$ belongs to
both of the above subcases, and that in each case $\nu=(21)$. This
case corresponds to the multiplicity $2$ appearing in the well
known expansion $s_{321/21}=s_3+2s_{21}+s_{111}$.

\noindent{\bf T2}.~~For $\nu=(xyzw)$ figure $K2$ shows the
preliminary constraints on interior edge labels that arise from
fixing the boundary edge labels.

Since $s_{\lambda/\mu}$ is basic, there are two subcases to deal
with, as tabulated below:
$$
\begin{array}{lclclcl}
{\bf T2}{\rm(i)}&&a>b>c\geq d+1\geq e>1
&&\nu=(a-1,b-1,c-e+1,d-e+1)\,;&&\cr
{\bf T2}{\rm(ii)}&&a>b>c\geq e\geq d+1>1
&&\nu=(a-1,b+d-e,c-e+1,0)\,.&&\cr
\end{array}
$$

For these $\nu$ we are able to complete the labelling of the
interior edges of two LR $3$-subhives of the required pair of LR
$4$-hives, as illustrated in the following figures:

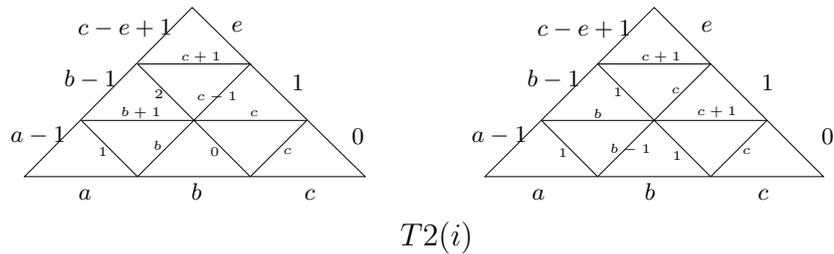
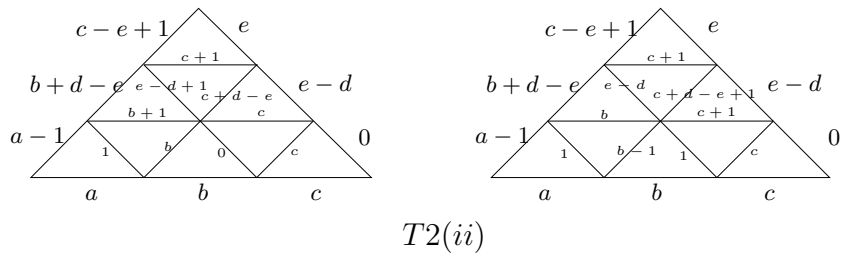
\begin{center}
\input{abcd-ee1.TpX}
\end{center}
\vspace{-10mm}
\begin{center}
\input{abcd-ee2.TpX}
\end{center}
\vspace{-5mm}

The edge labellings of these pairs of LR $3$-hives serve to complete
the interior edge labelling of the corresponding pairs of LR
$4$-hives in which they are embedded.
The existence of precisely two
LR-hives corresponding to the given $\nu$, shows that $s_{abcd/e^2}$
is not multiplicity-free.
Once again it might be noted that the case
$a=4$, $b=3$, $c=2$, $d=1$  and $e=2$ belongs to both subcases, and
that in each case $\nu=(321)$. This corresponds to the multiplicity
$2$ appearing in the expansion
$s_{4321/2^2}=s_{42}+s_{41^2}+s_{3^2}+2s_{321}+s_{31^3}+s_{2^3}+s_{2^21^2}$.

\noindent{\bf T3}.~~Setting $\nu=(wxyzuv)$, figure $K3$ illustrates
the impact of the specification of boundary edge labels on the
interior edges.

There are just two subcases to deal with,
and in each of these we consider $\nu=(w,x,y,z,u,v)$ as tabulated below:
$$
\begin{array}{lll}
\!\!{\bf T3}{\rm(i)}&\!\!a-1>b>c+1\geq d>2&\nu=(a-1,a-2,b-1,b-d+1,c-d+2,c-d+1)\,;\cr
\!\!{\bf T3}{\rm(ii)}&\!\!a-1>b>d\geq c+1>2&\nu=(a-1,a+c-d-1,b+c-d,b-d+1,1,0)\,.\cr
\end{array}
$$

In each subcase, we are then able to complete the labelling of the
interior edges of two LR $3$-subhives, as illustrated in the
following figures, where $f=d-c+1$ in the pair of subhives $T3(ii)$.

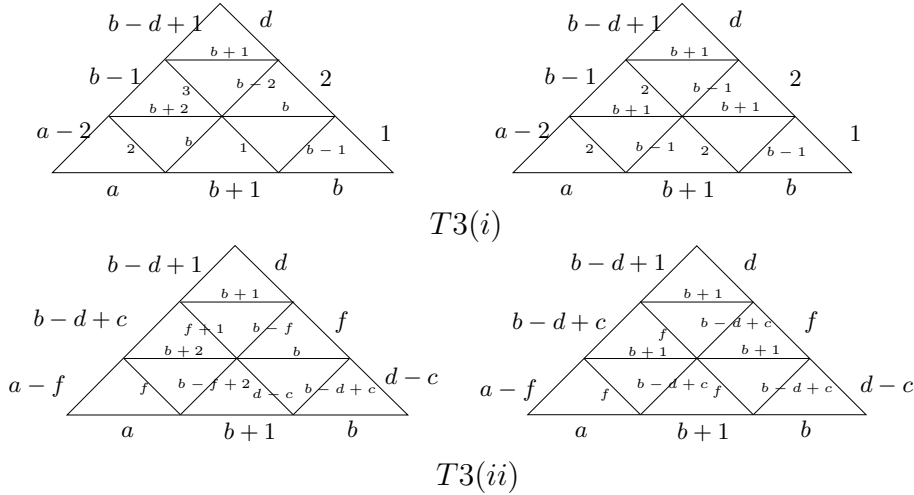
\begin{center}
\input{aabbcc-ddd1.TpX}
\end{center}
\vspace{-15mm}
\begin{center}
\input{aabbcc-ddd2.TpX}
\end{center}
\vspace{-8mm}

The edge labellings of these pairs of LR $3$-hives
serve to complete the interior edge labelling of the corresponding
pairs of LR $6$-hives in which they are embedded.
The existence of
precisely two LR-hives corresponding to the given $\nu$, then
suffices to show that $s_{a^2b^2c^2/d^3}$ is not multiplicity-free.
It might be noted that the two subcases coincide when $a=6$, $b=4$,
$c=2$ and $d=3$, in which case $\nu=(54321)$. This corresponds to
the multiplicity $2$ occuring in the decomposition
$s_{6^24^22^2/3^3}=s_{6531} + s_{6521^2} + s_{6432} + s_{6431^2} +
s_{642^2 1} + s_{6421^3} + s_{63^3}
 + s_{63^2 21} + s_{632^2 1^2} + s_{5^2 41} + s_{5^2 32} + s_{5^2 31^2} + s_{5^2 2^2 1}
 + s_{54^2 2} + s_{54^2 1^2} + s_{543^2} + 2s_{54321} + s_{5431^3} + s_{542^3} + s_{542^2 1^2}
 + s_{53^3 1} + s_{53^2 2^2} + s_{53^2 21^2} + s_{532^3 1} + s_{4^3 21} + s_{4^3 1^3}
 + s_{4^2 3^2 1} + s_{4^2 32^2} + s_{4^2 321^2} + s_{43^3 2} + s_{43^2 2^2 1}
$.

This completes the proof of Lemma~\ref{lem-multiplicity-main}. \qed

\noindent{\bf Note}:\ \ It should be pointed out that the conditions
on $a, b, c, d, e$ as stated in Lemma \ref{lem-multiplicity-main}
arise naturally. In $K1$, to avoid being
multiplicity-free, part (iii) of Lemma~\ref{lem-hexagon} implies
that $a>b>c$ and $d>e>0$, while part (ii) implies $c>0$. The
remaining conditions $a>d$ and $b>e$ of case {\bf T1} arise from the
fact that~\eqref{eq-hive-condition-2}
and~\eqref{eq-triangle-condition} imply $a\geq d+z\geq d$ and $b\geq
e$ with an interior edge label of the hexagon fixed to be $a$ or $b$
if either $a=d$ or $b=e$, respectively. In $K2$, part (iii) of
Lemma~\ref{lem-hexagon} implies that we require $a>b>c$, $c>d$ and
$e>p=e-d+w>0$ so that $e>1$ and $d>w$ with $w\geq0$. The remaining
condition $b>e$ of case {\bf T2} arises
because~\eqref{eq-hive-condition-2}
and~\eqref{eq-triangle-condition} imply $b\geq e+w\geq e$ and
setting $b=e$ would fix an interior edge label of the hexagon.
Finally in $K3$, to avoid being multiplicity-free,
Lemma~\ref{lem-hexagon} gives $a>p>b$, $b>q>c$ and
$d>r=d-c+u>s=d-c+v$, so that $a>b+1$, $b>c+1$ and $c>u>v$, with the
last implying $c>1$. In addition,~\eqref{eq-triangle-condition}
gives $b=q+s$, so that $b>q$ implies $s>0$. Then our condition
$d>r>s$ gives $d>2$, as required. The condition $q>c$ further
implies $b>d+v$, which yields the final condition $b>d$ of case
{\bf T3}.

\begin{lem}
The skew Schur functions $s_{\lambda/\mu}$ are not
multiplicity-free in each of the following cases:
\begin{description}
  \item[U1]~{\rm(i)} $\lambda=(a^2bc)$, $\mu=(de)$ with $a>b>c>0$, $d>e>0$
  and $a>d$,\\
  {\rm(ii)} $\lambda=(abc^2)$, $\mu=(d^2e)$ with $a>b>c>0$, $d>e>0$,
  $b>d$ and $c>e$;
  \item[U2]~{\rm(i)} $\lambda=(a^2bcd)$, $\mu=(e^2)$ with $a>b>c>d>0$
  and $a>e+1>2$,\\
  {\rm(ii)} $\lambda=(abcd^2)$, $\mu=(e^3)$ with $a>b>c>d>0$, $c>e>1$ and $d>1$;
  \item[U3]~{\rm(i)} $\lambda=(a^3b^2c^2)$ and $\mu=(d^3)$ with
  $a>b+1$, $b>c+1$, $c>1$ and $a>d+2>4$,\\
  {\rm(ii)} $\lambda=(a^2b^2c^3)$ and $\mu=(d^4)$ with
  $a>b+1$, $b>c+1$, $c>2$ and $b>d+1>3$.
\end{description}
\label{cor-multiplicity-main}
\end{lem}

\pf Once again we note that under the stated conditions
$s_{\lambda/\mu}$ is necessarily row-basic, but may not be basic.
However, in each case we can obtain $F^{\hat\lambda/\hat\mu}$
from $F^{\lambda/\mu}$ by the deletion of empty columns. This deletion procedure
is such that the pair $\hat\lambda$ and $\hat\mu$ necessarily
belong to the same case, {\bf U1}(i)--{\bf U3}(ii),
as the original pair $\lambda$ and $\mu$.
Since $s_{\lambda/\mu}=s_{\hat\lambda/\hat\mu}$, it follows
once again, that without loss of generality, we can confine attention
to those $s_{\lambda/\mu}$ that are basic. We consider each such case in turn.

{\bf U1}(i).~~Since $s_{\lambda/\mu}$ is basic, we have $b\geq e$.
If $b>e$ then the pair $\sigma=(abc)$ and $\tau=(de)$ are such that
$s_{\sigma/\tau}$ is row-basic. It follows from case {\bf T1} of
Lemma \ref{lem-multiplicity-main} that there exists at least one
$\rho$ such that $c^{\sigma}_{\tau\rho}\geq 2$. Then by the second
part of \eqref{eq-add-col-row} $c^{\lambda}_{\mu\nu}\geq
c^{\sigma}_{\tau\rho}\geq 2$, with $\nu=\rho\cup\{a\}$.

If $b=e$, we begin with $\sigma=(a-1,b,c)$ and $\tau=(d-1,b-1)$.
Case {\bf T1} of Lemma \ref{lem-multiplicity-main} applies to this
pair, so there exists at least one $\rho$ such that
$c^{\sigma}_{\tau\rho}\geq 2$. Then by the second part of
\eqref{eq-add-col-row}
$c^{\sigma^\prime}_{\tau^{\prime}\rho^{\prime}}\geq
c^{\sigma}_{\tau\rho}\geq 2$ with $\sigma^\prime=\sigma\cup\{a-1\}$,
$\tau^{\prime}=\tau\cup\{0\}$ and $\rho^{\prime}=\rho\cup\{a-1\}$.
At last by the first part of \eqref{eq-add-col-row}, we have
$c^{\lambda}_{\mu\nu}\geq
c^{\sigma^\prime}_{\tau^{\prime}\rho^{\prime}}\geq 2$ with
$\nu=\rho^{\prime}$.

{\bf U1}(ii).~~Since $s_{\lambda/\mu}$ is basic, we have $c\geq d$.
With $\sigma=(abc)$ and $\tau=(de)$ it follows from case {\bf T1} of
Lemma \ref{lem-multiplicity-main} that there exists at least one
$\rho$ such that $c^{\sigma}_{\tau\rho}\geq 2$. Then by the second
part of \eqref{eq-add-col-row} $c^{\lambda}_{\mu\nu}\geq
c^{\sigma}_{\tau\rho}\geq 2$, with $\nu=\rho\cup\{c-d\}$.

{\bf U2}(i).~~Since $s_{\lambda/\mu}$ is basic, we have $b\geq e$.
If $b>e$, let $\sigma=(abcd)$ and $\tau=(e^2)$. Then by case {\bf
T2} of Lemma \ref{lem-multiplicity-main} there exists at least one
$\rho$ such that $c^{\sigma}_{\tau\rho}\geq 2$. The second part of
\eqref{eq-add-col-row} then implies $c^{\lambda}_{\mu\nu}\geq
c^{\sigma}_{\tau\rho}\geq 2$ with $\nu=\rho\cup\{a\}$.

If $b=e$, let $\sigma=(a-1,b,c,d)$ and $\tau=((b-1)^2)$. By case
{\bf T2} of Lemma \ref{lem-multiplicity-main} there exists at least
one $\rho$ such that $c^{\sigma}_{\tau\rho}\geq 2$. Then by the
second part of \eqref{eq-add-col-row}
$c^{\sigma^\prime}_{\tau^{\prime}\rho^{\prime}}\geq
c^{\sigma}_{\tau\rho}\geq 2$ with $\sigma^\prime=\sigma\cup\{a-1\}$,
$\tau^\prime=\tau\cup\{0\}$ and $\rho^\prime=\rho\cup\{a-1\}$. At
last by the first part of \eqref{eq-add-col-row}, we have
$c^{\lambda}_{\mu\nu}\geq
c^{\sigma^\prime}_{\tau^{\prime}\rho^{\prime}}\geq 2$ with
$\nu=\rho^{\prime}$.

{\bf U2}(ii).~~Since $s_{\lambda/\mu}$ is basic, we have $d\geq e$.
Let $\sigma=(abcd)$ and $\tau=(e^2)$, so that by case {\bf T2} of
Lemma \ref{lem-multiplicity-main} there exists at least one $\rho$
such that $c^{\sigma}_{\tau\rho}\geq 2$. Then by the second part of
\eqref{eq-add-col-row} $c^{\lambda}_{\mu\nu}\geq
c^{\sigma}_{\tau\rho}\geq 2$ with $\nu=\rho\cup\{d-e\}$.

{\bf U3}(i).~~Since $s_{\lambda/\mu}$ is basic, we have $b\geq d$.
If $b>d$, let $\sigma=(a^2b^2c^2)$ and $\tau=(d^3)$. By case {\bf
T3} of Lemma \ref{lem-multiplicity-main} there exists at least one
$\rho$ such that $c^{\sigma}_{\tau\rho}\geq 2$. Then by the second
part of \eqref{eq-add-col-row} $c^{\lambda}_{\mu\nu}\geq
c^{\sigma}_{\tau\rho}\geq 2$ with $\nu=\rho\cup\{a\}$.

If $b=d$, let $\sigma=((a-1)^2b^2c^2)$ and $\tau=((b-1)^3)$. By case
{\bf T3} of Lemma \ref{lem-multiplicity-main} there exists at least
one $\rho$ such that $c^{\sigma}_{\tau\rho}\geq 2$. Then by the
second part of \eqref{eq-add-col-row}
$c^{\sigma^\prime}_{\tau^{\prime}\rho^{\prime}}\geq
c^{\sigma}_{\tau\rho}\geq 2$ with $\sigma^\prime=\sigma\cup\{a-1\}$,
$\tau^\prime=\tau\cup\{0\}$ and $\rho^\prime=\rho\cup\{a-1\}$. At
last by the first part of \eqref{eq-add-col-row}, we have
$c^{\lambda}_{\mu\nu}\geq
c^{\sigma^\prime}_{\tau^{\prime}\rho^{\prime}}\geq 2$ with
$\nu=\rho^{\prime}$.

{\bf U3}(ii).~~Since $s_{\lambda/\mu}$ is basic, we have $c\geq d$.
Let $\sigma=(a^2b^2c^2)$ and $\tau=(d^3)$. By case {\bf T3} of
Lemma~\ref{lem-multiplicity-main} there exists at least one $\rho$
such that $c^{\sigma}_{\tau\rho}\geq 2$. Then by the second part of
\eqref{eq-add-col-row} $c^{\lambda}_{\mu\nu}\geq
c^{\sigma}_{\tau\rho}\geq 2$ with $\nu=\rho\cup\{c-d\}$. \qed

The significance of these results is that it allows us to prove the
main theorem.

\noindent{\bf Proof of Theorem~\ref{thm-main}}\ \ Let $\lambda$ and
$\mu$ be such that $s_{\lambda/\mu}$ is basic, with $\lambda$ and
$\mu$ having $s$ and $t$ distinct non-zero parts, respectively, with
$s\geq1$ and $t\geq0$. Then for $s>0$ and $t>0$ we let
$\lambda=(a_1^{p_1},a_2^{p_2},\ldots,a_s^{p_s})$ and
$\mu=(b_1^{q_1},b_2^{q_2},\ldots,b_t^{q_t})$, with
$a_1>a_2>\cdots>a_s>0$, $b_1>b_2>\cdots>b_t>0$,
$\ell(\lambda)=p_1+p_2\cdots+p_s=n$ and
$\ell(\mu)=q_1+q_2+\cdots+q_t<n$, where $p_i,q_j>0$ for
$i=1,2,\ldots,s$ and $j=1,2,\dots,t$.

First we recall that the results of Section~\ref{sec-main-result}
imply that $s_{\lambda/\mu}$ is multiplicity-free in each of the
cases {\bf S0}--{\bf S8}. Then we
consider all possible values of $s$ and $t$ in turn.

In the following, we select the parts of $\sigma$ and $\tau$
from those of $\lambda$ and $\mu$, respectively, in such
a way that if $\sigma_i=\lambda_j$ then $\tau_i=\mu_j$.
Since $s_{\lambda/\mu}$ is basic, this guarantees that
$s_{\sigma/\tau}$ is at least row-basic. If
$s_{\sigma/\tau}$ is not multiplicity-free, there exists at
least one $\rho$ such that $c^{\sigma}_{\tau\rho}\geq 2$. Setting
$\lambda=\sigma\cup\zeta$ and $\mu=\tau\cup\xi$, we now
define $\eta$ to be the partition such that each
$\eta_k=\zeta_l-\xi_l$ for some $l$, and let
$\nu=\rho\cup\eta$. It then follows
from the repeated application of the second part
of \eqref{eq-add-col-row} with $a=\zeta_l$, $b=\xi_l$ and
$c=\eta_k$, that $c^\lambda_{\mu\nu}\geq c^\sigma_{\tau\rho}\geq2$.
Thus $s_{\lambda/\mu}$ is not multiplicity-free.

{\it The case $t=0$ is covered for all $s$ by ${\bf S0}$, and is
multiplicity-free}. If $s\geq3$ and $t\geq2$. We select the $\sigma$
and $\tau$ according to the relations between the various $p_i$ and
$q_i$. Three situations may arise: (i) If $\ell(\mu)\leq p_1$ then
we can always select $\sigma=(a^2bc)$ and $\tau=(de)$ which is
covered by case {\bf U1}(i) of Lemma~\ref{cor-multiplicity-main};
(ii) If $p_1<\ell(\mu)$ and $q_1<n-p_s$ then we can select
$\sigma=(abc)$ and $\tau=(de)$ which is covered by case {\bf T1} of
Lemma~\ref{lem-multiplicity-main}; (iii) Finally, if $p_1<\ell(\mu)$
and $q_1\geq n-p_s$ then we can select $\sigma=(abc^2)$ and
$\tau=(d^2e)$ which is covered by case {\bf U1}(ii) of
Lemma~\ref{cor-multiplicity-main}.

If $s\geq4$ the case $t\geq2$ has already been dealt with. We can
therefore take $t=1$ so that $\mu=(e^q)$ with $1\leq q<n$. If $q=1$
or $e=1$ or $e=a-1$ or $q=n-1$ then $s_{\lambda/\mu}$ is
multiplicity-free since the situation is covered by case {\bf S1},
which was dealt with in Section~\ref{sec-main-result}. For $a>e+1>2$
and $1<q<n-1$, we need only consider the following subcases: (i) If
$2\leq q\leq p_1$ we can select $\sigma=(a^2bcd)$ and $\tau=(e^2)$
which is covered by case {\bf U2}(i) of
Lemma~\ref{cor-multiplicity-main}; (ii) If $p_1<q<\ell(\lambda)-p_s$
we can select $\sigma=(abcd)$ and $\tau=(e^2)$ which is covered by
case {\bf T2} of Lemma~\ref{lem-multiplicity-main}; (iii) If
$\ell(\lambda)-p_s\leq q\leq n-2$ we can select $\sigma=(abcd^2)$
and $\tau=(e^3)$ which is covered by case {\bf U2}(ii) of
Lemma~\ref{cor-multiplicity-main}.

The case $s=1$ is covered for all $t$ by ${\bf S0^{\pi}}$, and is
multiplicity-free. Similarly, the case $s=2$ and $t=1$ is covered by
{\bf S4}, and is multiplicity-free. By virtue of the rotation
symmetry~\eqref{eq-skew-rotate}, the case $s=2$ and $t=2$ is
identical with that of $s=3$ and $t=1$, which is still to be fully
covered. On the other hand the cases $s=2$ and $t\geq3$ are
identical with their images under rotation, for which $s\geq4$ and
$t=1$, that have just been covered.

This just leaves the case $s=3$ and $t=1$. This is dealt with by
noting that if $\lambda=(a^{p_1}b^{p_2}c^{p_3})$ and $\mu=(d^q)$ are
such that they are not covered by one or other of the
multiplicity-free cases listed under {\bf S1}, {\bf S2} or {\bf S3},
that is $p_i\geq2$ for $i=1,2,3$, $a>b+1$, $b>c+1$, $c>1$, $d>2$,
$a>d+2$ and $3\leq q\leq n-3$, then they are covered by one of the
following three subcases: (i) If $3\leq q\leq p_1$ we select
$\sigma=(a^3b^2c^2)$ and $\tau=(d^3)$ which is covered by case {\bf
U3}(i) of Lemma~\ref{cor-multiplicity-main}; (ii) If $p_1<q<p_1+p_2$
we can select $\sigma=(a^2b^2c^2)$ and $\tau=(d^3)$ which is covered
by case {\bf T3} of Lemma~\ref{lem-multiplicity-main}; (iii) If
$p_1+p_2\leq q\leq n-3$ we can select $\sigma=(a^2b^2c^3)$ and
$\tau=(d^4)$ which is covered by case {\bf U3}(ii) of
Lemma~\ref{cor-multiplicity-main}.

This completes the proof of Theorem~\ref{thm-main}. \qed

\section{Final remarks}
\label{sec-final-remarks}

We have shown that the hive model is well suited to the derivation
of the two Theorems~\ref{thm-stem} and~\ref{thm-main} on
multiplicity-free Schur function products and skew Schur functions,
respectively. The use of LR-hives has allowed a direct proof that
all the cases enumerated in both theorems are indeed
multiplicity-free. In addition it has enabled us to demonstrate that
the breakdown of multiplicity-freeness always has a common origin,
in the sense that it can be traced back to the existence of a vertex
in the relevant LR-hives that is surrounded by an elementary
hexagon, none of whose interior edge labels is fixed either by the
criteria of Lemma~\ref{lem-hexagon} or any other means.

The proof offered here of the skew Schur function theorem, unlike
that of Gutschwager~\cite{guts}, is quite independent of
Stembridge's product of Schur functions theorem. In fact, since the
hive model proof has covered simultaneously both connected and
disconnected cases, it is possible to recover from
Theorem~\ref{thm-main} not only Theorem~\ref{thm-stem} but
also the following:

\begin{cor}
\label{cor-skew-products}
The product $s_\theta\,s_\phi$ of any two basic skew Schur functions
is multiplicity-free if and only if one or more of the following
is true:
\begin{description}
\item[V1]~~$\theta$ is a one-line rectangle and $\phi$ or $\phi^\pi$
is a partition (or vice versa);
\item[V2]~~$\theta$ is a two-line rectangle and $\phi$ or $\phi^\pi$
is a fat hook (or vice versa);
\item[V3]~~$\theta$ is a rectangle and $\phi$ or $\phi^\pi$ is a near-rectangle (or vice versa);
\item[V4]~~$\theta$ and $\phi$ are rectangles.
\label{description-skewprodlist}
\end{description}
\end{cor}

\pf\ Thanks to~\eqref{eq-disconnected} every product
$s_\theta\,s_\phi$ of two basic skew Schur functions can be
expressed as a single basic skew Schur function $s_{\lambda/\mu}$
where $F^{\lambda/\mu}$ is constructed, as in
Example~\ref{exmp-disconnected}, by joining $F^\theta$ and $F^\phi$
corner to corner. Then one applies Theorem~\ref{thm-main} to all
those cases for which $F^{\lambda/\mu}$ is of the required
disconnected form. The cases {\bf R0} are always
connected, while each multiplicity-free disconnected case of
Theorem~\ref{thm-main} gives rise through the identity
$s_{\lambda/\mu}=s_\theta\,s_\phi$ to a corresponding case in this
corollary, and {\it vice versa}, as follows:
${\bf R1}\leftrightarrow{\bf V1}$;
${\bf R2}\leftrightarrow{\bf V2}$;
${\bf R3}\leftrightarrow{\bf V3}$
and
${\bf R4}\leftrightarrow{\bf V4}$.
\iffalse
${\bf
S2}\leftrightarrow{\bf V1}$; ${\bf S3}\backslash{\bf
S8}\leftrightarrow{\bf V1}^\pi$; ${\bf S4}\backslash{\bf
S6}\leftrightarrow{\bf V2}$; ${\bf S5}\backslash{\bf
S7}\leftrightarrow{\bf V2}^\pi$; ${\bf S6}\backslash{\bf
S2}\leftrightarrow{\bf V3}$; ${\bf S7}\backslash{\bf
S3}\leftrightarrow{\bf V3}^\pi$; ${\bf S8}\leftrightarrow{\bf V4}$.
where $\,^\pi$ indicates that it is necessary to $\pi$-rotate either
$\theta$ or $\phi$. Using~\eqref{eq-skew-rotate} where necessary,
\fi
It follows that every multiplicity-free skew Schur function product is
of one or other of the types {\bf V1}--{\bf V4}. All other cases are
not multiplicity-free. \qed

By exactly the same argument, one can recover Theorem~\ref{thm-stem}
as a further corollary by restricting attention to those cases for
which both $\theta$ and $\phi$ are partitions.
%%% without recourse to any rotation through $\pi$.
The correspondence between the multiplicity-free disconnected cases
of Theorem~\ref{thm-main} of the required form and those of
Theorem~\ref{thm-stem} is given by: ${\bf R1}\leftrightarrow{\bf
P1}$; ${\bf R2}\leftrightarrow{\bf P2}$; ${\bf
R3}\leftrightarrow{\bf P3}$ and ${\bf R4}\leftrightarrow{\bf P4}$.
All other cases are not multiplicity-free.

\iffalse
It might, finally, be noted that the extension of these ideas to the
determination of multiplicity-free products of Schubert classes has
already been carried out by Gutschwager~\cite{guts}, who reproduced
the result of Thomas and Yong~\cite{tho-yon}. This can be done
equally well from our statement of Theorem~\ref{thm-main} which is
well suited to the case of the basic, but possibly disconnected skew
diagrams that may arise.
\fi

\noindent{\bf Acknowledgement}\ \ This article is based on work
carried out by two of us, DQJD and RLT, in partial fulfilment of
study for PhD degrees. RCK is pleased to acknowledge the hospitality
extended to him by Professor Bill Chen and his colleagues during
several visits to the Center for Combinatorics at Nankai University.
We are also indebted to Professor Bessenrodt for drawing to our
attention the work of her student Gutschwager.

\end{document}